\documentclass[12pt]{amsart}
\usepackage[english]{babel}
\usepackage{color}
\usepackage[all]{xy}
\usepackage{faktor}
    \usepackage{amsfonts}
    \usepackage{amsmath,amstext,amssymb,amsthm,amsfonts}
    \usepackage{array}
    \usepackage{tikz}
    \usepackage{mathtools}
    \usetikzlibrary{chains}
    \usepackage{cancel}
\usepackage{color}
\usepackage{tikz-cd}

\newcommand{\HA}{\mathcal{H}}

\newcommand{\mfg}{\mathfrak{g}}
\newcommand{\mfk}{\mathfrak{k}}
\newcommand{\mfp}{\mathfrak{p}}
\newcommand{\mfm}{\mathfrak{m}}
\newcommand{\mfa}{\mathfrak{a}}

\newcommand{\mfl}{\mathfrak{l}}
\newcommand{\mfh}{\mathfrak{h}}
\newcommand{\mfq}{\mathfrak{q}}
\newcommand{\mft}{\mathfrak{t}}

\newcommand{\CC}{\mathbb{C}}
\newcommand{\HH}{\mathbb{H}}
\newcommand{\RR}{\mathbb{R}}
\newcommand{\ZZ}{\mathbb{Z}}
\newcommand{\NN}{\mathbb{N}}

\newcommand{\CP}[1]{\mathbb{C}P^{#1}}
\newcommand{\HP}[1]{\mathbb{H}P^{#1}}

\newtheorem{te}{Theorem}[section]
\newtheorem{pr}[te]{Proposition}
\newtheorem{co}[te]{Corollary}
\newtheorem{lm}[te]{Lemma}

\newtheorem{re}[te]{Remark}
\newtheorem{ex}[te]{Example}

\newcommand{\Tr}{\operatorname{Tr}}
\newcommand{\Span}{\operatorname{Span}}
\newcommand{\rk}{\operatorname{rk}}

\title[Laplace eigenfunctions and Borel-Weil Theorem]{Laplace eigenfunctions on Riemannian symmetric spaces and Borel-Weil Theorem}

\begin{document}

\author{Dimitar Grantcharov, Gueo Grantcharov and Camilo Montoya}
\date{\today}
\address{ D. Grantcharov: Department of Mathematics\\
University of Texas Arilington\\
 Arlington, TX 76019-0408\\
 USA }
 \email{grandim@uta.edu}
 \address{ G. Grantcharov and C. Montoya: Department of Mathematics and Statistics,  Florida International University\\
Miami, FL 33199\\
   USA}
\email{grantchg@fiu.edu, camimont@fiu.edu}

\maketitle

\begin{abstract}
We identify a geometric relation between  the Laplace-Beltrami spectra and eigenfunctions on compact Riemannian symmetric spaces and the Borel-Weil theory using ideas from symplectic geometry and geometric quantization. This is done by associating to each compact Riemannian symmetric space, via Marsden-Weinstein reduction, a generalized flag manifold which covers the space parametrizing all of its maximal totally geodesic tori. In the process we notice a direct relation between the Satake diagram of the symmetric space and the painted Dynkin diagram of its associated flag manifold. We consider in detail the examples of the classical simply-connected spaces of rank one and the space $SU(3)/SO(3)$. In the second part of the paper, with the aid of harmonic polynomials, we induce Laplace-Beltrami eigenfunctions on the symmetric space from holomorphic sections of the associated line bundle on the generalized flag manifold. In the examples we consider we show that our construction provides all of the eigenfunctions.\\



\end{abstract}

\section{Introduction}

%\textcolor{blue}{We have different notations for the scalar products in the different sections. Should we explain them in each section separately or not? This might be better once we use real complex and quaternionic notations. And just mention where they are fixed it in this section.}

There are two classical geometric interpretations of the representation theory of the compact Lie groups. On the one side is the Borel-Weil Theorem and its subsequent generalization to the Borel-Weil-Bott theory. In particular, every complex representation of a compact Lie group is realized on the space of holomorphic sections of some line bundle over a flag manifold. On the other side, the harmonic analysis on a Riemannian symmetric spaces provides irreducible representations of a compact simple Lie group:  the natural action of a transitive simple Lie group of isometries $G$ on the common eigenspaces of the (commutative) algebra of the invariant differential operators on the respective compact symmetric space $M=G/K$ is irreducible. This algebra contains the Laplace-Beltrami operator and is generated by $k$ generators, where $k=\rk(M)$ is the rank of $M$. One of the goals of this paper, of which some results were reported in \cite{GG}, is to propose a direct geometric relation between the two theories. The most explicit illustration of the relation is through geometric quantization of the geodesic flow in the case $\rk(M)=1$ and we briefly explain it first.

We assume for simplicity that such compact Riemannian symmetric space of rank one (CROSS for short) is simply-connected and irreducible, this leaves the well-known examples of the spheres, classical projective spaces (complex and quaternionic) and the Cayley plane. These are also the known simply-connected examples in dimension higher than two of Riemannian manifolds all of whose geodesics are closed. On a Riemannian manifold $(M, g)$  all of whose geodesics are closed there is a natural $S^1$-action on its tangent bundle $TM$ and the geodesic flow on the cotangent bundle $T^*M$ can be realized as solution to an $S^1$-invariant Hamiltonian system. For such systems,  under mild conditions, there is a moment map and a symplectic reduction process, called also Marsden-Weinstein reduction.  This reduction produces a reduced space $T^*M// S^1$ that can be identified with the space parametrizing all oriented geodesics and that is equipped with an induced symplectic form. The induced symplectic form depends on a level set of the corresponding moment map which  we call the energy level of the geodesic flow. In many examples the cotangent bundle has a "complex polarization" - a complex structure compatible with the symplectic form which becomes K\"ahler form. We use this form to apply a twisted version of Kostant-Souriau geometric quantization scheme (originally due to \cite{Cz, He}, )and assign a holomorphic line bundle with first Chern class given by the induced K\"ahler form with an added extra term. This new term is half of the first Chern class of the canonical bundle of the manifold $T^*M//S^1$. Sometimes this is called {\it prequantum bundle} and we use this terminology. The quantum condition is the integrality of that corrected form, while the analog of the Hilbert space of quantum observables is the space of holomorphic sections of the prequantum bundle. Our first general result is Theorem \ref{rk1} in Section 4 in which we show that the quantized energy levels of the geodesic flow on a simply connected rank one symmetric space are, up to a constant, equal to the eigenvalues of the Laplace-Beltrami operator on $M$, and the corresponding complexified eigenspaces are isomorphic to the
spaces of the holomorphic sections of the prequantum bundle over the reduced space which is a generalized flag manifold. Although we provide the representation theory background for a unified proof, we proceed with a case by case proof since it illustrates the explicit nature of the calculations and provides a basis for the next parts of the paper.

 In the next Sections we consider the case of general rank. We observe that we can substitute the space parametrizing all geodesics with the space of all maximal totally geodesic flat submanifolds, which are tori in this case. Just as in the rank one case we needed the oriented geodesics, in the higher rank case we need the universal cover of the space parametrizing the maximal flat tori. This space is again a generalized flag manifold and carries a natural "polarization" which could be used for the quantization - a K\"ahler complex structure. Since our aim is to underline the geometric approach through the Marsden-Weinstein reduction, we also need a  K\"ahler space with a (multi-dimensional) Hamiltonian that, after the symplectic reduction, will become the generalized flag manifold with an appropriate reduced symplectic form, a form which is  also integral and K\"ahler.
This is done in \cite{GG} via construction of a K\"ahler structure on some open subset of the manifold of all tangent spaces of the maximal totally geodesic flat submanifold. Then we prove the main result in Section 6 - Theorem \ref{main} which was announced in \cite{GG}. If the symmetric space has a maximal rank $\rk(M) = \rk(G)$, then the corresponding generalized flag manifold is actually the full flag manifold $G/T$, where $T$ is a maximal torus in $G$. From the Borel-Weil theorem follows that every irreducible representation of $G$ appears as a space of holomorphic sections of some line bundle over $G/T$. This corresponds to the fact that the symmetric spaces of maximal rank provide the largest variety of the irreducible representations of $G$ appearing as subspaces of the eigenspaces of the Laplace-Beltrami operator by the results in \cite{H2}. As an example we consider the space $SU(3)/SO(3)$  which is of rank two, and also the simplest example of Riemannian symmetric space of maximal rank. We note that a general correspondence similar to the one in Theorem \ref{main} has appeared in \cite{Gi} and \cite{H3} Chapter 6, but with a focus on various integral transforms.

In the remaining part of the paper we indicate a construction which relates the holomorphic sections of the prequantum bundle to the eigenfunctions of the corresponding eigenvalue on the symmetric space. We use the standard description of the holomorphic sections as holomorphic functions on the total space of the associated principal $\CC^*-$bundle with appropriate equivariant condition. Using the basic properties of the Laplace-Beltrami operator under Riemannian submersions, the relation between holomorphic and harmonic functions on a special non-K\"ahler manifolds, and an extension of the standard relation between harmonic polynomials and eigenfunctions on the spheres, in Theorem \ref{rel} we propose a method of description of the Laplace-Beltrami eigenfunctions out of the holomorphic sections on the associated  quantization space. We apply the method to complex and quaternionic projective spaces as well as the space $SU(3)/SO(3)$. In these examples we describe a spanning set of all eigenfunctions, which consists of algebraic functions. The cases of quaternionic spaces and $SU(3)/SO(3)$ are new and extend the known representations for the spheres and complex projective spaces \cite{Ta}. We expect that many other symmetric spaces will have similar complete description. Note that the known descriptions of the eigenfunctions so far are based on the integral geometry and various types of Radon transform which are implicit.

The structure of the paper is as follows: in the preliminary sections we collect the necessary background facts about symmetric spaces and flag manifolds and notice a simple connection between the Satake diagram of the symmetric space and the painted Dynkin diagram of the generalized flag manifold. We use it do describe the second cohomology group of the quantization space in terms of the Satake diagram of the initial symmetric space. Then in Section 4 we treat the rank one case and in Sections 5 and 6 we consider the spaces of any rank. The last two Sections - 7 and 8 provide the construction of harmonic polynomials from holomorphic sections of the prequantum line bundle over the corresponding generalized flag manifold.

\medskip
\noindent \textbf{Acknowledgements:} The first author is partly supported by  Simons Collaboration Grant 358245. He also would like to thank the Max Planck Institute in Bonn (where part of this work was completed) for the excellent working conditions. We are grateful to Francis Burstall for the stimulating discussion and valuable remarks.  We thank also Joe Wolf and Simon Salamon for their inetrest.

\section{Riemannian symmetric spaces and generalized flag manifolds}

Let $G$ be a compact semisimple Lie group and $K$ a Lie subgroup given by a fixed point set of an automorphic involution $\theta$. Then $M=G/K$ carries a Remannian metric which makes it a Riemannian symmetric space.  Denote by $\mfg$ and $\mfk$ the Lie algebras of $G$ and $K$. There is an eigenspace decomposition of $\theta$ $\mfg=\mfk + \mfp$ where we identify $\mfp$ with $T_oM$ the tangent space at the $o=eK$ of $M$. It is known that $\mfg_n = \mfk + i\mfp$ defines the non-compact dual Lie algebra of $\mfg$ with respect to $\theta$. Denote by $G_n$ the corresponding simply-connected Lie group. It is known that $[\mfk,\mfp]= \mfp$ and  $[\mfp,\mfp]=\mfk$. Denote by $\mfa$ the maximal abelian subalgebra in $\mfp$. Then dimension of $\mfa$ is called a {\it rank} of $M$. Denote also by $\mfg^\CC, \mfk^\CC$, $\mfp^\CC$ etc. the complexifications of $\mfg, \mfk,\mfp$ etc. respectively. Note that $\mfg^\CC$ is a complexification of $\mfg_n$ and $i\mfa $ is a maximal abelian in $i\mfp$.

The non-compact group $G_n$ admits an Iwasawa decomposition $G_n=KAN$ where $A$ is the simply-connected Lie group with algebra $i\mfa$ and $N$ is unipotent. There is a complex Iwasawa decomposition (see e.g. \cite{Cl}) given by  $G^\CC_0=K^\CC A^\CC N^\CC$ where $K^\CC,  A^\CC, N^\CC$, and $G^\CC$ are the complexifications of $K,A,N,G$ and $G^\CC_0$ is some Zarisky open and dense subset of $G^\CC$.

It is known that $K$ acts transitively on the set of all maximal abelian subalgebras in $\mfp$. Denote by $M$ the stabilizer of $\mfa$ in $K$ and by $\mfm$ its Lie algebra. It follows that it may be expressed as $\mfm= \{X\in\mfk : [X,\mfa]=0\}$ and since $\mfa$ is maximal, the centralizer of $\mfa$ in $\mfg$ is $\mfl=\{X\in\mfg : [X,\mfa]=0\}=\mfm + \mfa$. If $L$ is the corresponding subgroup in $G$, then it contains a maximal torus of $G$ (as a centralizer of an abelian subgroup) so the space $G/L$ is a generalized flag manifold (also referred to in the literature as a K\"ahlerian C-space and rational homogeneous manifold) and carries a natural complex structure, as well as a K\"ahler metric. Its geometric interpretation is elucidated in the following comment.
 \begin{re} From geometric viewpoint a maximal abelian subalgebra of $\mfp$ is tangent to maximal totally geodesic flat torus and every such torus is at a point $gK$ is tangent to a left translate of some $\mfa$ from $o=eK$ to $gK$. In particular the space parametrizing all maximal totally geodesic flat tori can be identified with $G/N_{\mfa}$ where $N_{\mfa}$ is the normalizer of a fixed $\mfa$ in $\mfg$. As shown in \cite{H2} the centralizer $L$ is a normal subgroup of $N_{\mfa}$ and the quotient $N_{\mfa}/L$ is a finite group which is a subgroup of the Weyl group of $\mfg$. In this terms the generalized flag manifold $G/L$ is the universal cover of the space parametrizing the set of all such tori. We call the space $G/L$ a {\bf quantization space} of $G/K$. It is closely related to the horospherical manifold in \cite{H3} and \cite{Gi}.
 \end{re}

%\vspace{1.469mm}
An important invariant related to the quantization space is the dimension of its second cohomology. We demonstrate here how this dimension could be identified in terms of the data provided by the Satake diagram associated to the noncompact dual Riemannian symmetric space of $G/K$. Denote by $\mfm_0=[\mfm,\mfm]$. Since $\mfm$ is compact, $\mfm_0$ is the semisimple part of $\mfm$. We can consider a maximal torus of $\mfg$ which is $\theta$ invariant and contains $\mfa$. Such torus is known to exists and since $\mfa$ is maximal abelian in $\mfp$, then this torus is a Cartan subalgebra and has the form $\mfh = \mft + \mfa$ with $\mft \subset\mfk$. Its complexification $\mfh^{\CC} = \mft^{\CC}+\mfa^{\CC}$ is a Cartan subalgebra of $\mfg^{\CC}$. Let $\Delta = \Delta_+ \cup \Delta_-$ be a root system with an ordering defining the positive and negative roots of $\mfg^{\CC}$.  There is a set $\Sigma$ of the so called restricted roots $\Delta \supset \Sigma = \Sigma(\mfg^\CC,\mfa^\CC)\subset (\mfa^\CC)^*$ and we can choose a basis $h_1,...,h_k,h_{k+1},...,h_n$ of $(\mfh^\CC)^*$, of basic roots, such that $h_1,...,h_k$ (after restricting them via the projection $\mfh^\CC\rightarrow\mfa^\CC$) are basis for $(\mfa^\CC)^*$. We continue to use the same notation $h_1,...,h_k$ for the restricted roots. After we choose an ordering of the basic roots, or equivalenlty a positive Weyl chamber, every element of $\Delta$ (respectively, $\Sigma$( is an integer linear combinations with all non-negative or all non-positive coefficients of $h_1,...,h_n$ (respectively, $h_1,...,h_k$). Similarly we can choose a positive Weyl chamber in the restricted roots.  Then $\mfm_0^{\CC}$ has a root decomposition with root spaces which are among the root spaces $\mfg_{\alpha}$ of $\mfg^{\CC}$ with respect to $\mfh^{\CC}$.

\section{Satake and painted Dynkin diagrams}

 The non-compact real form $\mfg_n$ of $\mfg^{\CC}$ which is dual of $\mfg$ have $\mfk$ as a maximal compact subalgebra. We want to describe the relation between the Satake diagram of $\mfg_n$ and the painted Dynkin diagram associated to the generalized flag manifold $G/L$ (recall that  $\mfl = \mfa + \mfm$). We first describe the Satake diagram of $\mfg_n$. For this we paint the simple roots defining the root spaces of $\mfm_0^{\CC}$ in black. The Satake diagram of $\mfg_n$ is then the Dynkin diagram of $\mfg^{\CC}$ with the black and white dots as described, but with an additional arrows between the white roots, when there is an involutive automorphism of $\mfg$ such that the difference between the initial and the endpoint of an arrow is a root in $\mft$ ( see \cite{OV} Sect. 4.4 or \cite{CS}  Sect. 2.3,  for the facts about Satake diagrams).

If we consider a parabolic subalgebra $\mfq$ of $\mfg^{\CC}$ which contains $\mfl^{\CC} = \mfm^{\CC} + \mfa^{\CC}$ and is minimal (i.e. does not contain another proper parabolic subalgebra), then it defines the complex structure of the generalized flag manifold $G/L$. Such algebra is not unique, but the different ones are related via elements of the Weyl group and they define different invariant complex structures on $G/L$. As described in \cite{BR}, $\mfq$ contains a Borel subalgebra and is described by a subset of simple roots of $\Delta_+$ and the complex structure is defined via the ordering.  The painted Dynkin diagram of $G/L$ in an analog of the Satake diagram  and is defined  - see  \cite{AlP, BFR} as a diagram with vertices corresponding to the semisimple part of $\mfl$ painted in black. Since $\mfl = \mfa+\mfm$, the semisimple (ss) part of $\mfl$ is precisely $\mfl_{ss} = \mfm_0$. The description leads to the following observation:
\begin{te}
Let $M = G/K$ is an irreducible Riemannian symmetric space and $G/L$ is its quantization space. If $\mfg_n$ is the non-compact dual form of $\mfg$ with respect to $K$, then the Satake diagram of $\mfg_n$ with the arrows deleted corresponds precisely to the painted Dynkin diagram defining the complex structure of $G/L$.
\end{te}

From here we obtain:

\begin{co}\label{H2Sat}
The group $H^2(G/L, \ZZ)$ has no torsion and is generated by the elements corresponding to the white vertices of the Satake diagram of $\mfg_n$.
\end{co}

{\it Proof:} The fact that $H^2$ has no torsion is well known. According to \cite{BH}, there is an isomorphism
${\mathcal Z}(\mfl^{\CC})^* \equiv H^2(G/L, \mathbb{R})$ sometimes called {\it transgression}, given by $\alpha \rightarrow \frac{i}{2\pi} d\alpha$. It has the property, that the fundamental weights in $\mathcal{Z}(\mfl^{\CC})^*$ correspond to elements of $H^2(G/L,
\ZZ) \in H^2(G/L,\RR)$. Then the first $k$ elements $w_1,...,w_k$  of the basis $w_1,...,w_n$ of the fundamental weights dual w.r.t Killing form on $\mfg^{\CC}$ to $h_1,...h_n$, define an integral basis of $H^2(G/L,\mathbb{Z})$. It is known \cite{AlP, BFR}, that this basis is generated  via transgression by simple roots corresponding to the white vertices in the painted Dynkin diagram describing the complex structure of $G/L$ which are in bijection with the basis of the center of $\mfl$. As explained above these are precisely the white vertices in the corresponding Satake diagram.{\it Q.E.D.}

\begin{re}
The results clarifies a statement in Lemma 12 of \cite{GG}.

\end{re}

In what follows we use the notations $h_1,...,h_n$ and $w_1,...,w_n$ for the simple roots and dual positive weights throughout the paper, as well as the notion of transgression from the proof above.

\medskip
\noindent
{\bf Example.} Consider the complex Grassmannian $SU(p+q)/S(U(p)\times U(q))$ with $p\leq q$. Then $s(u(p)+u(q)) = su(p)+su(q)+\RR$ is block-diagonally embedded in $su(p+q)$. If $su(p+q)$ is represented as the matrices of the form
$\left(\begin{array}{ll}
A & B \\
-\overline{B}^t & C
\end{array}\right)$
with $B$ a complex $p\times q$ matrix  and $A,C$ skew-Hermitian, then $\mfk = s(u(p)+u(q))$ corresponds to $B=0$, $\mfp$ corresponds to $A=C=0$ which provides the decomposition $\mfg=\mfk+\mfp$. As is known (cf. \cite{H1,H2,H3}) and could be checked directly, the space $\mfa$ could be identified with the matrices $\left(\begin{array}{ll}
0 & D \\
-D^t & 0
\end{array}\right)$, where $D$ is a real-valued upper-diagonal matrix $\left(\begin{array}{lllllll}
d_1 & 0 &...&0 & 0 & ...& 0\\
0& d_2 & ... &0 & 0 & ...& 0\\
.& .&...&.\\
0& 0&...& d_p & 0 & ...& 0\\
\end{array}\right)$ with a real diagonal matrix on the first p rows and p columns, only $0$ entries on the remaining $p-q$ columns. Now we can check directly that $\mfm$
is the spaces of the trace-free  matrices in $\mfk$ (i.e. $B=0$), where $A$ is diagonal matrix with imaginary entries, and $C$ = $\left(\begin{array}{ll}
E & 0\\
0 & F
\end{array}\right)$, with $E$ being a  $p\times p$ diagonal matrix with imaginary entries and $F$ a skew-Hermitian $(q-p) \times (q-p)$ matrix. From here we can see that the center of $\mfl = \mfm+\mfa$ consists of $\mfa$ as well as the matrices  $\left(\begin{array}{lll}
D_1& 0 & 0\\
0& D_1 & 0\\
0 & 0& \alpha Id
\end{array}\right)$, with $D_1$ - diagonal with imaginary entries and $\alpha = -2\Tr(D_1)$ - an imaginary complex number. In particular $\mfm_0 = su(q-p)$ is the space with skew-Hermitian matrix in the lower-right corner. Then  the generalized flag manifold is $G/L \cong SU(p+q)/S(T^{2p}\times U(q-p))$ and the dimension of $H^2(G/L,\ZZ)$ is $2p$. Note that the Cartan subalgebra $\mfh$ in this case is different from the standard one consisting of all diagonal imaginary matrices in $su(p+q)$. If $\{e_1,...,e_{p+q}\} $ is the canonical basis in $\CC^n$, then the algebra $\mfh$ is represented by the diagonal  matrices in the basis $\{f_i = e_i + ie_{i+p}, \overline{f_i}, e_j\}$ where $1\leq i\leq p$ and $p+1\leq j\leq p+q$.

%We need some more information on the roots and weights of Lie algebra $\mfg^c$ relative to $\mfk^c$ and $\mfa^c$. Choose a maximal abelian (Cartan) %subalgebra $\mfh^c$ of $\mfg^c$ which contains $\mfa^c$. Denote by $\mfh^c_k$ the intersection $\mfh^c\cap \mfk^c$, so that $\mfh^c=\mfh^c_k+\mfa^c$ %and let $\Delta\subset (\mfh^c)^*$ be the corresponding root system of $\mfg^c,\mfh^c$. There is a set of the so called restricted roots $\Delta %\supset \Sigma = \Sigma(\mfg^c,\mfa^c)\subset (\mfa^c)^*$ and we can choose a basis $h_1,...,h_k,h_{k+1},...,h_n$ of $(\mfh^c)^*$, of basic roots, %such that $h_1,...,h_k$ (after restricting them via the projection $\mfh^c\rightarrow\mfa^c$) are basis for $(\mfa^c)^*$. We continue to use the %same notation $h_1,...,h_k$ for the restricted roots. After we choose an ordering of the basic roots , or equivalenlty a positive Weyl chamber, %every element of $\Delta$ resp. $\Sigma$ is an integer linear combinations with all non-negative or all non-positive coefficients of $h_1,...,h_n$ ( %$h_1,...,h_k$ resp.). Similarly we can choose a positive Weyl chamber in the restricted roots.

%An important observation is:

%\begin{lm}
%The center $\mathcal{Z}(\mfl^c)$ coincides with $\mfa^c$

%\end{lm}

%For the proof just note that $\mfa^c$ is contained in the center by definition, and any other hypothetical element should be in $\mfh^c_k$. But $K$ %is maximal compact in $G_n$ and $\mfg^c$ is semisimple and $\mfp^c$ is $K$-irreducible  - so ??? -({\bf Da se pogledne - mislia che e prosto})....

\vspace{.2in}

\section{Symplectic reduction and quantization of the geodesic flow of the symmetric spaces of rank one}

We first recall Marsden-Weinstein (or symplectic) reduction and then the modified Konstant - Souriau geometric quantization scheme (twisted half of the canonical bundle - see  Czyz \cite{Cz} and Hess \cite{He}). In \cite{MT} the authors related the energy spectrum of the quantized geodesic flow on a sphere with the eigenvalues of the Laplace-Beltrami operator. In what follows we first notice a generalization of this result from a representation theory viewpoint to all compact rank-one Riemannian symmetric spaces (CROSSes). Then we provide a similar detailed computations for two examples -  $\CP n$ and $\HP n$. The exposition follows \cite{GG} and is used later to determine a generating set of eigenfunctions defined by harmonic polynomials in an ambient space.

For more details on the Marsden-Weinstein (or symplectic) reduction we refer the reader, for example, to \cite{AM}, here we provide a partial review of the results we need later. . If $(M,\omega)$ is a symplectic manifold and $H$ is a function on $M$, then the vector field $X_H$ defined as $dH(Y)=\omega(X_H,Y)$ is called Hamiltonian vector field. We will call $H$ a Hamiltonian function and the triple $(M,\omega,H)$ - a Hamiltonian system. If $G$ is a group of symplectomorphisms then under mild conditions there is a map $\mu: M\rightarrow \mathfrak{g}^*$, defined by $$d\mu(X) = i_{X}\omega,$$
where $\mathfrak{g}$ is the Lie algebra of $G$ and $X\in \mathfrak{g}$ is identified with the induced vector field on $M$. When such $\mu$ exists, the action is called {\it Hamiltonian} and  the space $N= \mu^{-1}(c)/G$ is called the Marsden-Weinstein reduction or the symplectic reduction, where $c$ is a fixed element of the adjoint action of $G$ on $\mathfrak{g}^*$. We denote $N$ by $M//G$. The space $M//G$ inherits a natural symplectic form $\omega_{red} $ such that $ i^*(\omega)=\pi^*(\omega_{red})$ where $i:\mu^{-1}(c)\rightarrow M$ is the inclusion and $\pi:\mu^{-1}(c)\rightarrow N=\mu^{-1}(c)/G$ is the natural projection. The following results will be used repeatedly in the paper (see e.g. \cite{AM} for the proof).

\begin{pr} \label{mw-reduction}
If $(N,\omega_{red})$ is the symplectic reduction of $(M,\omega)$ under the action of a Lie group $G$ and $H$ is a G-invariant function on $M$, then there is a unique function $H_{red}$ on $N$ such that $\pi^*(H_{red}) = i^*(H)$. Moreover the flow of the vector field $X_H$ preserves $\mu^{-1}(c)$ and projects on $N$ to the flow of the vector fields $X_{H_{red}}$. Moreover, if we have a second Hamiltonian action of a Lie group $G_1$ on $M$ which commutes with the action of $G$, then the level sets of its moment map $\mu_1$ are $G-$invariant and $\mu_1|_{\mu^{-1}(c)} = \pi^*(\overline{\mu_1})$ where $\overline{\mu_1}$ is the moment map associated to the action of $G_1$ on $N$.
\end{pr}

\begin{pr}
If $G$ is a compact group of isometries acting freely on the Riemannian manifold $(M,g)$ and $N= M/G$ is the orbit space, then for the canonical symplectic forms $\Omega_M, \Omega_N$ on $T^*M, T^*N$, respectively, we have $T^*N = T^*M//G$ with $\Omega_N$ being the reduced form from $\Omega_M$.
\end{pr}

%We'll use this facts in the next Section and see generalize them later in the paper.

%\section{Symmetric spaces of rank one and quantization of the geodesic flow}

 The geodesic flow on a Riemannian manifold is represented as a Hamiltonian flow on its cotangent bundle. The cotangent bundle of each Riemannian manifold $(M,g)$ has a canonical symplectic form given in local coordinates as $\Omega = \sum dx_i\wedge dy_i$ where $(x_1,...,x_n)$ are local coordinates of $M$ and $(x_1,...,x_n,y_1,...,y_n)$ are the associated local coordinates of $T^*M$. Then the function $H(x,v) = \frac{1}{2}g(v,v)$ for $x\in M$ and $v\in T^*_xM$ has a Hamiltonian vector field $X_H$ and its flow lines project on $M$ to give the geodesics. In particular, $i_{X_H}\Omega = dH$. If all the geodesics of $M$ are closed, then they define an $S^1$-action on $T^*M$ with orbits $(c(t), g(c'(t))$ for a geodesic $c(t)$ and the dual 1-form $g(c'(t))$ of its tangent vector $c'(t)$. Which means that the moment map $\mu(x,v)$  at $(x,v) \in T^*M, v\in T^*_xM$ for this $S^1$-action is precisely $\mu= H$. When we fix the level $c$ of the moment map, the points in the reduced space $H^{-1}(c)/S^1$ represent (oriented) geodesics on $M$ with tangent vectors of length $c$. Which explains that $Geod(M) = T^*M//S^1$ as sets where $Geod(M)$ is the set of oriented geodesics on $M$. We note that the reduced form $\Omega_c$ from the canonical form on $T^*M$ depends on the choice of the level set $\mu^{-1}(c)$ for the moment map of the action $\mu$ (which is called the {\it energy of the geodesic flow}).

Now recall some facts about the quantization scheme of Konstant and Souriau with the amends of Czyz and Hess \cite{Cz, He}. Let $X$ be a compact K\"ahler manifold with K\"ahler form $\lambda$. We say that the holomorphic line bundle $L$ is a {\it quantum line bundle} if its first Chern class satisfies $$c_1(L)=\frac{1}{2\pi}[\lambda] - \frac{1}{2}c_1(X).$$
Thus $X$ will be quantizable if and only if $\frac{1}{2\pi}[\lambda] - \frac{1}{2}c_1(X)\in H^2(X,\ZZ)$ . The corresponding quantum Hilbert space is the (finite dimensional) linear space $H^0(X, \mathcal{O} (L))$. We want to apply the scheme to the space of geodesics of a Riemannian manifold all of whose geodesics are closed.

 Main examples of such manifolds are the compact Riemannian symmetric spaces of rank one (CROSS for short). Recall that for a compact irreducible Riemannian symmetric space $G/K$ with simple Lie group $G$ we associated a quantization space $G/L$ which covers the space parametrizing all maximal totally geodesic flat tori in $G/K$. The space $G/L$ is a generalized flag manifold, so smooth projective variety and from the description of its second cohomology we know that its Picard group contains the center of $\mfl = \mfa+\mfm$. In particular it contains the fundamental weights $w_1,...,w_k$ which correspond to the restricted roots for $\mfa$.  Now denote by $\mathcal{L} = \mathcal{L}_{i_1,...,i_k}$ the holomorphic line bundle on $G/L$ determined by $w=i_1w_1+...+ i_kw_k$, where $i_j\geq 0$. By Bott vanishing the higher cohomology of $\mathcal{L}$ are zero. The space $H^0(G/L, \mathcal{O}(\mathcal{L}))$ is a (unitary) representation of $G$ with highest weight $w$. The Borel-Weil theorem shows that the representation is irreducible if $w$ is dominant, and corresponds to the (unique) irreducible representation with highest weight $w$ \cite{Serre}.

On the other side,  the general theory for the Laplace spectrum on symmetric spaces (\cite{CW, Ta}) tells us that the eigenvalues are given by $\lambda = ||\rho((\mfa^\CC)^*) + w||^2-||\rho((\mfa^\CC)^*)||^2$ where $w$ is as before and $\rho((\mfa^\CC)^*)$ is the half sum of positive restricted roots of $\mfa^\CC$. When the center of $\mfl$ is $\mfa$,  $\rho$ represents one half of the first Chern class of $G/L$, so $\rho((\mfa^\CC)^*)+w$ is the first Chern class of $\mathcal{L}\otimes K^{\frac{1}{2}}$.

%\textcolor{blue}{ Maj e viarno slednoto: The corresponding dimensions are also equal due to the Borel-Weil-Bott. This is the representation-theoretical explanation of the geometric %quantization results considered later.}

In this Section we focus on the case of Riemannian symmetric spaces of rank one, since the correspondence in this case is most studied and related to the classical quantization of the geodesic flow. In the next sections we'll generalize the scheme to the symmetric spaces of higher rank.
When the rank of $M$ is one,  the space of the restricted roots $\Sigma$ is 1-dimensional as is the Weyl chamber in it. The set of fundamental weights in it is (see \cite{H1}):

$$
\Lambda^+ = \left \{\lambda \in \mfa^\CC | \frac{\langle \lambda,\psi \rangle}{\langle  \psi,\psi \rangle} \in \ZZ^+, \textrm{for}\hspace{.02in}\textrm{all}\hspace{.02in} \psi\in \Sigma \right \}
$$
 and in the rank one case is generated by a single element $\theta$ . The considerations above give the following result, which we will generalize in the next Section.

 \begin{te}\label{rk1} let $M = G/K$ be an irreducible simply-connected  compact Riemannian symmetric space of rank one (CROSS). Then up to re-scaling of the metric on $M$ the following are true:

 i) Under the transgression the reduced symplectic form $\Omega_c$ on $Geod(M)=G/L = T^*M//S^1$ corresponds to $\pi\sqrt{2c}\theta$ and with the choice of the positive Weyl chamber and complex structure as above, $c_1(G/L)$ corresponds to $N_M\theta$ for a positive integer $N_M$.

 ii) The quantum condition on $(Geod(M),\Omega_c)$ (i.e. $[\Omega_c]\in H^2(G/L, \ZZ)$) provides the following energy spectrum: $c_k = 1/2(N_M+2k)^2$.

 iii) The spectrum of the (semi)-Laplacian $\frac{1}{2}\Delta_M$ on $M$ is given by $\lambda_k = ||k\theta+\rho(\mfa^\CC)||^2-||\rho(\mfa^\CC)||^2$ and $c_k =  ||k\theta+\rho(\mfa^\CC)||^2$ where $\rho(\mfa^c)$ is the half-sum of the positive restricted roots of $\mfa^\CC$.

 iv) The multiplicities of $c_k$ and $\lambda_k$ coincide with the dimension of the (finite - dimensional) representation $L(k\theta)$ of $\mfg$ with highest weight $k\theta$ relative to $(\mfh,\Delta)$. Moreover the representation $L(k\theta)$ is isomorphic to both the (complex) eigenspace $\mathcal{L}^2(M)^{\lambda_k}$ of $\Delta_M$ corresponding to $\lambda_k$ and the quantization space $H^0(Geod(M),\mathcal{O}({\mathcal L}_k))$.

 \end{te}

{\it Proof:} The spaces in the Theorem are classified and are $S^n, \CP n , \HP n, Ca P^2$. In the two examples below we give a proof in case of $\CP n$ and $\HP n$. The case of $S^n$ is considered in \cite{MT}.
In \cite{MT} the space of oriented geodesics of $M = S^n$ is explicitly identified with the complex quadric in $\CP n$ via the Marsden-Weinstein reduction. It was noted that the energy levels of the moment map that satisfy a quantization condition coincide, up to an additive constant, with the eigenvalues of the Laplace-Beltrami operator and the their multiplicity are the same as the (complex) dimension of the holomorphic sections of the corresponding quantum bundle $L(k\theta)$ ( see also \cite{Ra} for related results).

Finally, consider the case  $CaP^2 = F_4/Spin(9)$. Since $M$ has rank one, then the reduction identifies the level set of the moment map with a spheric bundle over $M$. From \cite{Fu} Proposition 3.3 follows that it is diffeomorphic to $F_4/Spin(7)$. This gives an identification of the quantization space with $F_4/Spin(7)\times S^1$. Its painted Dynkin diagram from \cite{BFR}, Table 4 and Corollary \ref{H2Sat} follows that $H^2(F_4/Spin(7)\times S^1 , \ZZ) = \ZZ$, so the reduced form $\Omega_c$ is proportional to the generator. Since the generator corresponds to $\theta$ under the transgression, and the proportionality constant depends on $c$, and $i)$ follows. Then $ii)$ follows by the quantization condition, $iii)$ and $iv)$ by combining the results from \cite{CW} and Borel-Weil Theorem. {\it Q.E.D.}

We note that the simply-connected requirement could be lifted and similar statement could be stated for $\RR P^n$.

\begin{re}
The re-scaling factor mentioned in the Theorem could be different for the different spaces. It is known that the eigenvalues of the Laplace-Beltrami operator for $S^n$ and $\CP n$ are $k(n+k-1)$ and $4k(n+k)$ in the round metric on $S^n$ and the Fubini-Studi metric on $\CP n$ respectively. So the two metrics are rescaled differently, one by a factor 4 times the other.

\end{re}
We continue with the explicit calculations of the two classical simply-connected projective spaces.

\subsection{Complex projective space} \label{subsec-cpn}
We'll use and extend here the results of \cite{FT}. Throughout this subsection for complex vectors $z,w$  we denote by $\langle z, w \rangle = \mbox{Re}\, \sum z_i\overline{w}_i$ their hermitian scalar product and by $z.w = \sum z_iw_i$ the complex scalar product so $||z||=\sqrt{\langle z, z\rangle}$. For a point $[z]=[z_0,z_1,...,z_n]$ in the complex projective space $\CP n$, we identify the holomorphic cotangent space
$$
T^*_{[u]}\CP n \cong \{ ([u],v)\in \{[u]\}\times \CC^{n+1}|  u.\overline{v}=0\}$$
where we used the Fubini-Study metric to identify the tangent and cotangent bundles. To achieve a global description of the cotangent bundle, we use the Hopf map $\pi: S^{2n+1}\rightarrow \CP n$ which is induced by the standard action of $S^1$ on $S^{2n+1}$. This map is defined by $u \mapsto [u]$, where $u\in \RR^{2n+2}=\CC^{n+1}$ with $||u|| = 1$. After identifying the tangent and cotangent bundles of the sphere via the canonical metric, we can identify the  cotangent bundle as
$$ T^*S^{2n+1} = \{ (u,v)\in \CC^{n+1}\times\CC^{n+1}| \hspace{.1in}||u||=1, \langle u,v \rangle = 0\}$$
Then the $S^1$-action $\rho$ for the Hopf projection $\pi$ extends to $T^*S^{2n+1}$ as $$\rho(e^{i\theta})(u,v)=(e^{i\theta} u,e^{i\theta}v)$$ This action preserves the canonical symplectic form on $T^*S^{2n+1}$, which is given by $i^*\mbox{Re}\, (du\wedge d\overline{v})$. The moment map for the action $\rho$ can be used to show the following theorem. This theorem is first proven in \cite{FT}, but for reader's convenience a short proof is presented. We consider the cotangent bundle with its zero section deleted $T^*_0\CP n$ (and $T^*_0S^{2n+1}$) in order to avoid the singularity issues since they are irrelevant in the paper.
\begin{lm} \label{diffeo}
 The space $T^*_0\CP n$ is diffeomorphic to  both $X_C$ and $\widetilde{X_C}$ where
 $$X_C \cong \{ [u,v]\; |\hspace{.1in} ||u||=1, u.\overline{v} = 0, v\neq 0\}$$
with $[u,v]$ representing the class of $(u,v)$ under $(u,v)\sim (e^{i\theta}u,e^{i\theta}v)$ and
$$\widetilde{X_C} \cong  \{ [[u,v]]\; |\hspace{.1in} \langle u,u \rangle= \langle v,v \rangle \neq 0,  u.v =0\}$$
with $[[u,v]]$ defined by the relation $(u,v)\sim (e^{i\theta}u,e^{-i\theta}v)$. Moreover $T^*_0\CP n$ is biholomorphic to $\widetilde{X_C}$ when it is identified with $T^*_0S^{2n+1}//S^1$ and the reduced complex structure.
\end{lm}
{\it Proof:} It is well-known that under the action $\rho$,  $T^*\CP n = T^* S^{2n+1}//S^1$. The moment map $\Phi$ associated to the action $\rho$ is simply $\Phi(u,v) = \mbox{Im}\langle u,v\rangle$. Hence, $T^*_0S^{2n+1}//S^1 = \Phi^{-1}(\mu)/S^1$, for a generic $\mu\in \RR=iu(1)$, is identified with $X_C$ which gives the diffeomorphism $T^*_0\CP m \cong X_C$.  The diffeomorphism between $X_C$ and $\tilde{X_C}$ is given by the formulas:
$$\tilde{u}_k=\frac{1}{\sqrt{2}}(||v||u_k+ iv_k),$$

$$\tilde{v}_k=\frac{1}{\sqrt{2}}(\overline{v}_k-i ||v|| \overline{u}_k).$$

The biholomorphism follows from the fact that the reduction is K\"ahler, when  we consider the canonical form on $T^*S^{2n+1}$ as a K\"ahler form for the complex structure induced from the embedding in $\CC^{2n+2}$ as in \cite{Ra} for example.

{\it Q.E.D.}

By lemma \ref{mw-reduction}, if  a Lie group $G$ of isometries acts on $M$, this action  induces a Hamiltonian action on $T^*M$ and the reduced space $T^*M//G$ becomes a (reduced) Hamiltonian system. Whenever $T^*M//G=T^*N$ for some Riemannian manifold $N$ then the solutions of the new system is precisely  the geodesic flow on $N$. In the particular case of $T^*S^{2n+1}$ we obtain the following.

\begin{pr}
The canonical symplectic form $\Omega_C$ on $T^*\CP n \cong X_C$ is $$\Omega_C = \frac{1}{2}(du\wedge d\overline{v}+d\overline{u}\wedge dv)$$ and the Hamiltonian system $\HA_{\CP n} = (X_c, \Omega_C, H_C= \frac{ ||v||^2}{2})$ induces the geodesic flow on $\CP n$. The system is equivalent to $(\widetilde{X}_C, \widetilde{\Omega}_C, \widetilde{H}_C) $ in view of the diffeomorphism in Lemma \ref{diffeo}.

\end{pr}

%{\color{red} Why $H_C= \frac{1}{ ||v||^2}$???}

Since the orbits of $\HA_{\CP n}$ correspond precisely to the geodesics of $\CP n$, we first identify the space parametrizing the geodesics. For this we first consider the geodesic flow on the sphere $S^{2n+1}$. Since all of the geodesics on the sphere are closed, the flow of $X_H$ in the cotangent space has also only closed trajectories. They define an $S^1$-action which is given by $(u,v)\rightarrow(e^{i\theta}u,e^{-i\theta}v)$. This action commutes with the action inducing the Hopf projection and is Hamiltonian. So it defines an action on $T^*\CP n$ which has orbits - the flow lines of the Hamiltonian vector field defining the geodesics on $\CP n$. We can identify a geodesic $c(t)$ in $\CP n$ with the line $(c(t),c'(t))$ in $T\CP n \simeq T^*\CP n$ when $t$ is a parameter such that $c'$ has constant norm. From here we see that the space parametrizing the geodesics can be identified with the Marsden-Weisntein quotient. Let $N_c = \widetilde{H}_C^{-1}(c)/S^1$ be the
reduced space. To identify $N_c$ with a flag manifold, we use the Hamiltonian system  $(\widetilde{X}_C, \widetilde{\Omega}_C, \widetilde{H_C}) $. Let

$$
\begin{array}{lll}
\mathbb{F}& = & \{ ([z],[w])\in \CP n\times \CP n :  z.w=0\} \\
& = & \{([z],[w])\in \CP n\times \CP n :  \langle z,z\rangle=\langle w,w\rangle, z.w =0\}
\end{array}
$$

One can see that $\mathbb{F}$ is biholomorphic to the (1,2)-flag in $\CC^{n+1}$ with homogeneous representation $\mathbb{F} = U(n+1)/U(1)\times U(1)\times U(n-1)$. Denote by $p_1$ and $p_2$ the two projections on the corresponding factors of $\CP n\times \CP n$. Let $\alpha$ be the generator (the Fubini-Study form) of $H^2(\CP n, \ZZ)$. Then $\omega_1 = p_1^*\alpha$ and $\omega_2 = p_2^*\alpha$ are generators of $H^2(\mathbb{F}, \ZZ)$. With this notation we have the following.

\begin{pr}
If $c\neq 0$ then the reduced manifold $N_c$ is biholomorphic to the flag $\mathbb{F}$ and the  reduced K\"ahler form is $\widetilde{\omega}_c = \pi\sqrt{2c}(\omega_1+\omega_2)$.

\end{pr}

{\it Proof:} The $S^1$-action of the geodesic flow on $T^* \CP n$ is induced from the one on $T^*S^{2n+1}$. Hence this action is:
$$
\lambda(z,w) = (\lambda z,\lambda w),$$
for $(z,w) \in \widetilde{H}_C^{-1}(c)$.  For the sphere $S^{2n+1}_R$ of radius $R$ the Hopf projection fits in the diagram $\xymatrix{\CC^{n+1}   & S^{2n+1}_R  \ar[l]_{i} \ar[r]^{h}  & \CP n}$
with $h^*\alpha = \frac{1}{\pi R^2}i^*\Omega$ (see \cite{MT}). If $\tilde{\pi}_c$ is the projection $\widetilde{H}_C^{-1}(c) \rightarrow N_c = \mathbb{F}$ then we have the following commutative diagram:
$$
\xymatrix{ \widetilde{H}_C^{-1}(c)  \ar[d]^{} \ar[r]^{\tilde{\pi}_c}  &  N_c \ar[d] ^{\tilde{i}_{c}}\\
S^{2n+1} \times S^{2n+1}  \ar[r]^{ h \times h}  &  {\mathbb C}P^n \times {\mathbb C}P^n },
$$
where the vertical arrows correspond to the natural embeddings. Therefore,
$$
\begin{array}{lll}
\tilde{\pi}^*_c(\sqrt{2c}\pi(\omega_1+\omega_2))& = &\pi\sqrt{2c}\frac{\sqrt{-1}}{2\pi}(\frac{dz\wedge d\overline{z}}{ ||z||^2}+ \frac{dw\wedge d\overline{w}}{ || w ||^2}) \\
& = & \frac{1}{\sqrt{2c}}\frac{\sqrt{-1}}{2}(dz\wedge d{\overline{z}} + dw\wedge d\overline{w})\\
& = & \frac{1}{2}(du\wedge d\overline{v} + d\overline{u}\wedge dv)\\
& = & {\tilde i}_{c}^*(\tilde{\Omega}_c)\\
\end{array}
$$

%\vspace{2in}

%{\bf Popravki ot 17 iuli}
%\vspace{.2in}

In the above calculation we used that $\widetilde{H}_C(z,w)=c$, so $||z||^2=||w||^2=2c$. We want to use the  modified Kostant - Souriau scheme to "quantize" the geodesic flow of $\CP n$.

\begin{pr}
We have $c_1(\mathbb{F})= n(\omega_1+\omega_2)$

\end{pr}

{\it Proof:} We apply the adjunction formula for a hypersurface of degree (1,1) in $\CP n\times \CP n$ to obtrain
$$
\begin{array}{lll}
c_1(\mathbb{F})&= &-(c_1(K_{\CP n\times \CP n}|_{\mathbb{F}})+c_1([\mathbb{F}]|_{\mathbb{F}}))\\
& = & c_1(\CP n\times \CP n)|_{\mathbb{F}}-c_1([\mathbb{F}]|_{\mathbb{F}})\\
& = & (n+1)(\omega_1+\omega_2) - (\omega_1+\omega_2)\\
& = & n(\omega_1+\omega_2)
\end{array}
$$

Q.E.D.

\begin{te} \label{th-cp}
The energy spectrum of the geodesic flow on $\CP n$ is: $$E_k = \frac{1}{2}(n+2k)^2, k\in \NN,$$ with corresponding multiplicities $$m_k = \binom{n+k}{k}^2 - \binom{n+k-1}{k-1}^2.
$$
\end{te}

{\it Proof:} For the exact cohomology sequence: $$H^1(\mathbb{F},\mathcal{O}) \rightarrow H^1(\mathbb{F},\mathcal{O}^*)\rightarrow H^2(\mathbb{F},\ZZ)\rightarrow H^2(\mathbb{F},\mathcal{O})$$ and the identities $ H^(\mathbb{F},\mathcal{O})=H^2(\mathbb{F},\mathcal{O})=0$ follows that:
$$c_1: H^1(\mathbb{F},\mathcal{O}^*)\xrightarrow{\cong} H^2(\mathbb{F},\ZZ)\cong \ZZ\oplus\ZZ.$$
Therefore every holomorphic line bundle $L$ on $\mathbb{F}$ is equivalent to $L_{k_1,k_2}=k_1\pi_1^*(H) + k_2\pi_2^*(H)$, where $H$ is the hyperplane section on $\CP n$.

The quantum condition on $c$ is : $$ \frac{1}{2\pi}[\omega_c]-\frac{1}{2}c_1(\mathbb{F})=c_1(L_{k_1,k_2})$$
which implies $$\frac{\sqrt{2c}}{2}-\frac{n}{2}=k$$ where $k=k_1=k_2$ is a positive integer. In particular $$c=\frac{1}{2}(2k+n)^2$$

To count the multiplicities (i.e. $\dim H^0(\mathbb{F},\mathcal{O}(L)$) we consider the exact sequence of sheaves:
\begin{equation}\label{S1}
 0\rightarrow \mathcal{O}_{\CP n\times \CP n}(L_{k-1,k-1}\otimes L_{1,1})\xrightarrow{\alpha} \mathcal{O}_{\CP n\times\CP n}(L_{k,k})\xrightarrow{r} \mathcal{O}|_{\mathbb{F}}(L_{k,k})\rightarrow 0,
\end{equation}
where $\alpha$ is the multiplication of sections of $L_{k,k}$ by the polynomial $\sum_0^nz_iw_i$ which defines $\mathbb{F}$ in $\CP n\times\CP n$ and $r$ is the restriction.
The corresponding exact cohomology sequence gives:
\[
0\rightarrow H^0(\CP n\times\CP n, \mathcal{O}(L_{k-1,k-1}))\rightarrow H^0(\CP n\times\CP n, \mathcal{O}(L_{k,k}))
\]
\[
\rightarrow H^0(\mathbb{F},\mathcal{O}(L_{k,k}))\rightarrow H^1(\CP n\times\CP n, \mathcal{O}(L_{k-1,k-1}))=0
\]
where the last term is zero by the Kodaira vanishing theorem. Thus we have:

$$
\begin{array}{lll}
m_k& = & \dim(H^0(\mathbb{F}, \mathcal{O}(L_{k,k}))\\
& = &\dim(H^0(\CP n\times\CP n, \mathcal{O}(L_{k,k})) - \dim(H^0(\CP n\times\CP n, \mathcal{O}(L_{k-1,k-1})))\\
& = &  \binom{n+k}{k}^2 - \binom{n+k-1}{k-1}^2.
\end{array}
$$
Q.E.D.

\subsection{Quaternionic projective space}

We first note  that the results in this subsection were independently obtained in \cite{El} and some of them appear in \cite{FT}. For readers convenience, in this section we use a slightly different notations to distinguish between real complex and quaternionic scalar products. In particular we use  $\langle x,y \rangle_\RR$ $\langle x,y \rangle_\CC$ and $\langle x,y \rangle_\HH$ for $ \overline{x}.y = \sum\overline{x_i}y_i$ when $x_i, y_i $ are in $\RR, \CC$, and  $\HH$ respectively. The corresponding norms arising from their rea parts are denoted by $||.||_\RR, ||.||_\CC, ||.||_\HH$ respectively. The geodesic flow on $\HP n$ can be described in a similar way as the one for $\CP n$ but with the aid of the quaternionic Hopf map.  For that we use three equivalent representations of $T^*S^{4n+3}$:
$$
\begin{array}{lll}
T^*S^{4n+3}& = & \{ (x,y)\in \RR^{4n+3}\times \RR^{4n+3}: ||x||_\RR=1, \langle x,y \rangle_\RR=0\}\\
& = &\{ (u,v)\in \CC^{2n+2}\times \CC^{2n+2}: ||u||_\CC=1, \mbox{Re} \langle u,v \rangle_\CC=0\}\\
& = & \{ (p,q)\in \HH^{n+1}\times \HH^{n+1}: ||p||_\HH=1, \langle p,q \rangle_\RR=0\}\\
\end{array},
$$
where $p_k:=u_{2k}+u_{2k+1}j, q_k:=v_{2k}+v_{2k+1}j$ and $\langle p,q\rangle_\RR = \mbox{Re} \langle p,q\rangle_\HH= \mbox{Re} \sum\overline{p}_kq_k$. The quaternionic Hopf map in this case is $\chi: S^{4n+3}\rightarrow \HP n, \, p\rightarrow[p]$ where $[p]=[p_0,p_1,...,p_n]$ is the class of $p$ for the relation $p\sim\sigma p, \sigma\in Sp(1)$. The next lemma is again from \cite{FT}.
\begin{lm}
The cotangent space $T^*\HP n$ is diffeomorphic to both $X_H$ and $\widetilde{X}_H$ defined as follows:
$$
X_H:= \{ \lfloor p,q\rfloor\in \HH^{n+1}\times \HH^{n+1}: ||p||_\HH=1, \langle p,q \rangle_\HH=0\},
$$
$$
\widetilde{X}_H:=\{\lfloor z,w\rfloor\in\CC^{2n+2}\times\CC^{2n+2}; ||z||_\CC=||w||_\CC, \langle z,w \rangle_\CC=0, I(z,w)=0 \}
$$
where $I(z,w)=z_0w_1-z_1w_0+...+z_{2n}w_{2n+1}-z_{2n+1}w_{2n}$ and $\lfloor p,q\rfloor$ and $\lfloor z,w \rfloor$ denote the equivalence classes of $(p,q)$ and $(z,w)$ under $(p,q)\sim (\sigma p,\sigma q)$ and $(z,w)\sim (z,w)g$ for $\sigma \in Sp(1)$ and $g\in SU(2)\cong Sp(1)$.
\end{lm}
{\it Proof:} Consider the action of $SU(2)$ on $S^{4n+3}$ defined by
\begin{equation}\label{psi}
\Psi_g(p,q):=(p,q)g, g\in SU(2).
\end{equation}
This action has a moment map $G: T^*S^{4n+3}\rightarrow su(2)^*$,  given by the formulas $G(p,q)=(A(p,q),B(p,q),C(p,q))$, where $$\langle p,q\rangle _\HH=\mbox{Re}(\langle p,q\rangle_\HH)+A(p,q)i+B(p,q)j+C(p,q)k,$$  and  the imaginary quaternions are identified with $su(2)^*$. Hence, $T^*S^{4n+3}//SU(2) = X_H\cong T^*\HP n$.

To prove that $X_H$ and $\widetilde{X}_H$ are diffeomorphic, consider  the map $t_H: X_H\rightarrow \tilde{X}_H, (z,w)=t_H(p,q)$, where
$$z_{2k}:=\frac{1}{\sqrt{2}}(||v||_{\mathbb C}u_{2k}+\sqrt{-1}v_{2k})$$
$$z_{2k+1}:=\frac{1}{\sqrt{2}}(-||v||_{\mathbb C}\overline{u}_{2k+1}-\sqrt{-1}\overline{v}_{2k+1})$$
$$w_{2k}:=\frac{1}{\sqrt{2}}(v_{2k+1}-\sqrt{-1}||v||_{\mathbb C}u_{2k+1})$$
$$w_{2k+1}:=\frac{1}{\sqrt{2}}(\overline{v}_{2k+1}-\sqrt{-1}||v||_{\mathbb C}\overline{u}_{2k+1})$$
Q.E.D.

The action $\Psi$ defined in (\ref{psi}) commutes with the geodesic flow of $S^{4n+3}$. Recall the diffeomorphism $t_H: X_H \to \widetilde{X}_H$ defined at the end of the last proof. Like in the previous subsection, we have the following.
\begin{pr}\label{pr1}
Let $\Omega_H = \Omega_{T^*\HP n}$ be the canonical symplectic form on $T^*\HP n$. Then $$\Omega_H = \frac{1}{2}(du\wedge d\overline{v}+d\overline{u}\wedge dv)$$
Moreover the geodesic flow of $\HP n$ is the flow of the equivalent Hamiltonian systems $$(X_H,\Omega_H,G_H)\cong(\widetilde{X}_H,\widetilde{\Omega}_H, \widetilde{G}_H)$$ where $G_H = \displaystyle \frac{||q||_\HH^2}{2} = \frac{||v||_\CC^2}{2}$,
%{\color{red} Maybe better write $G_H = \displaystyle \frac{||v||_C^2}{2}$?}
 $\widetilde{\Omega}_H=t_H^*\Omega_H$ and $\widetilde{G}_H = t_H^*({G}_H)$.
\end{pr}

Next we compute  the energy spectrum of the geodesic flow on $\HP n$ in a similar way as in the case of $\CP n$. We consider again the reduced space $\mathbb{O}_c = T^*\HP n//S^1 = \tilde{G}^{-1}(c)/S^1$ with the induced symplectic form $\omega_c$ obtained from $\tilde{i}^*_c\Omega_H=\tilde{\pi}_c^*\omega_c$, where $\tilde{i}_c:\widetilde{G}^{-1}(c)\rightarrow T^*\HP n$ and $\tilde{\pi}_c:\widetilde{G}^{-1}(c) \rightarrow \mathbb{O}_c$. Denote by $\mathbb{F}_{is}$ the isotropic Grassmann manifold
$$
\begin{array}{lll}
\mathbb{F}_{is}& = &\{\Lambda\in Gr_2(\CC^{2n+2}); I|_{\Lambda}=0\}\\
& = &\{[[z,w]] \in \CC^{n+1}\times\CC^{n+1} :  ||z||_\CC=||w||_\CC=1, \langle z,w \rangle_\CC=I(z,w)=0\}
\end{array}
$$
where $[[z,w]]$ is representative of $(z,w)$ for $(z,w)\cong (\lambda z,\lambda w)g$, $\lambda\in S^1, g\in SU(2)$
or equivalently $(z,w)\cong (z,w)g, g\in U(2)$. Alternatively, $\mathbb{F}_{is}$ is a hyperplane in $Gr_2(\CC^{2n+2})$:
$$
\mathbb{F}_{is}\cong \{(\lambda_{ij})\in Gr_2(\CC^{2n+2})\, : \lambda_{01}+\lambda_{23}+...+\lambda_{2n+1, 2n+2}=0\},
$$
where $(\lambda_{ij}) $ are the Pl\"{u}cker coordinates on $Gr_2(\CC^{2n+2})$, as well as  a  homogeneous space: $\mathbb{F}_{is}\cong Sp(n+1)/U(2)Sp(n-1)$.
\begin{pr}\label{pr2}
If $c\neq 0$ then the reduced space $\mathbb{O}_c$ is isomorphic to $\mathbb{F}_{is}$ equipped with the K\"ahler form $\widetilde{\omega}_c =\pi\sqrt{2c}\omega$, where $\omega$ is the restriction of the canonical K\"ahler form on $Gr_2(\CC^{2n+2})$ which generates $H^2(Gr_2(\CC^{2n+2}), \ZZ)$.
\end{pr}

{\it Proof:} The $S^1$ action of the geodesic flow on $\widetilde{G}_H^{-1}(c)\subset T^*\HP n \cong \widetilde{X}_H$ is:
$$\lambda \lfloor z,w\rfloor =\lfloor \lambda z, \lambda w\rfloor.$$
which commutes with the action of $Sp(1)\cong SU(2)$  defining the  quaternionic Hopf fibration. Now from $\widetilde{G}_H(z,w)=c$ we have $||z||_\CC^2=||w||_\CC^2=2c$. If $\lambda_{ij} = z_iw_j-z_jw_i$ are the Pl\"{u}cker coordinates on $Gr_2(\CC^{2n+2})$ then:
$$
\begin{array}{lll}
\tilde{\pi}^*_c(\pi\sqrt{2c}\omega)&=&\pi\sqrt{2c}\frac{\sqrt{-1}}{2\pi}\frac{d\lambda_{ij}\wedge d\overline{\lambda}_{ij}}{\sum_{i,j}||\lambda_{ij}||^2}\\
&=&\frac{1}{\sqrt{2c}}\frac{\sqrt{-1}}{2}(dz\wedge d\overline{z}+dw\wedge d\overline{w})\\
&=&\tilde{i}^*_c(\widetilde{\Omega}_H)
\end{array}
$$
Q.E.D.

\begin{pr}\label{pr3}
We have $c_1(\mathbb{F}_{is})= (2n+1)\omega$.
\end{pr}
{\it Proof:} We note that $c_1(Gr_2(\CC^{2n+2})|_{\mathbb{F}_{is}}=(2n+2)\omega$ and then proceed with the adjunction formula as in Proposition 2.4 using the fact that $\mathbb{F}_{is}$ is a hypersurface in $Gr_2(\CC^{2n+2})$.  Q.E.D.

\begin{te} \label{th-hp}
The energy spectrum of the geodesic flow on $\HP n$ is
$$E_k=\frac{1}{2}(2n+1+2k)^2, k\in \mathbb{N}$$
with corresponding multiplicities:
$$
m_k = \frac{2n+2k+1}{(k+1)(2n+1)}\binom{2n+k}{k}\binom{2n+k-1}{k}.
$$

\end{te}

{\it Proof:} We only  sketch the proof since it is similar to the  $\CP n$ case. We have $c_1:H^1(\mathbb{F}_{is},\mathcal{O}^*) \rightarrow H^2(\mathbb{F}_{is},\mathbb{Z}) = \mathbb{Z}$. Therefore all holomorphic line bundles on $\mathbb{F}_{is}$ which arise from the quantization are $L_k:=S^{\otimes k}$, where $S=\iota^*([H])$ and $\iota$ is the inclusion $\iota: \mathbb{F}_{is}\rightarrow Gr_2(\CC^{2N+2})$). Hence,
$$
\frac{\sqrt{2c}}{2} - \frac{2n+1}{2} = k.
$$
The dimension can be calculated via the Weyl dimension formula (see for example \cite{CW}).
{\it Q.E.D.}
\vspace{.2in}

This Theorem finishes the case by case proof of Theorem \ref{rk1}.

\begin{re}Existence of  a symplectic form on the cotangent bundle  doesn't have a direct analog to use in the case of higher-rank symmetric spaces. However not all simple and simply connected Lie group acts transitively on a CROSS. And the representation theory suggests that the correspondence could be extended to the higher dimensional case. In the next Sections we present one possible extension of the correspondence to the higher rank spaces based on geometry of toric bundles over flag manifolds.
\end{re}

%The most related to this is the paper by S. Gindikin arXiv:math/0501022 . There from the complex Iwasawa decomposition one can see that $G/L = G^c/M^cA^cN^c$ as a flag and $\Theta = G^c/M^cN^c$ is a complex space which is important for a CAuchy-Radon transform. The space $\Theta$ is precisely the space parametrizing all tangent spaces to all totally geodesic tori in $G/K$. Its "symplectic quotient" by the $k$-tori (max totally geodesics) will produce $G/L$. This is because in general there is a GIT quotient by $A^c$, which is the same as symplectic reduction for the Kahler case. It is known the $\Theta$ is Stein - i.e. has an exact Kahler form. But how to choose a particular one which will work - I am not sure. Some invariance should be involved. The best case scenario - it will depend on $k$ parameters, and their "quantum" condition should be - they are proportional to positive integers.

\section{Symplectic geometry of complex torus bundles}

Let $T^n=S^1\times S^1\times...\times S^1$ be the (real) n-dimensional torus. Then its tangent bundle is trivial and there is a well-known  identification $T(T^n) \equiv T(S^1)\times T(S^1)\times ...\times T(S^1) \equiv (\CC^*)^n\equiv (T^n)^\CC$ which is the complex n-dimensional torus. In particular it is an open and dense subset of $\CC^n$ and has an induced  complex structure and K\"ahler metric. If $z_k=r_ke^{i\theta_k}$ are the coordinates in $(\CC^*)^n$ then the K\"ahler form can be written as $\omega=\sum d(r_k^2)\wedge d\theta_k = d(\sum r_k^2d\theta_k)$. One can see that the $T^n$ action on $(\CC^*)^n$ $(z_1,...,z_n)\rightarrow (e^{i\alpha_1}z_1,...,e^{i\alpha_n}z_n)$ is Hamiltonian with moment map $\mu(z_1,...,z_n)=(r_1^2,...,r_n^2)=(|z_1|^2,...,|z_n|^2)$. Now we want to extend it to torus bundles:

\begin{te}\label{toricsymplectic} Let $\pi: P\rightarrow M$ be a principal $T^n$-bundle over a K\"ahler manifold $M$ with characteristic classes of type $(1,1)$. Let $P^\CC = P\times_{T^n}(\CC^*)^n$ be the associated complex torus bundle with the standard right action of $T^n$ on $(\CC^*)^n$. then $P^\CC$ is open and dense subset of the vertical tangent bundle $\mathcal{V}$ of $P$ and carries a natural complex structure and compatible symplectic (pseudo-K\"ahler)  form $\omega$. Moreover the $T^n$ action on $P^\CC$ is Hamiltonian and the Marsden-Weinstein reduction $P^\CC//T^n$ is diffeomoerphic to $M$ for a generic level set of the corresponding moment map.
\end{te}

{\it Proof:} A principal torus bundle is determined, up to an isomorphism, by its characteristic classes on the base. Then we have a closed and integral $(1,1)$-forms on $M$, $\omega_1,...,\omega_n$ and a connection 1-forms $\theta_1,..,\theta_n$ on $P$, such that $d\theta_k=\pi^*(\omega_k)$. The projection map $z:P\times (\CC^*)^n\rightarrow (\CC^*)^n$ defines functions $r_k^2=|z_k|^2$ which are $T^n$-invariant and descend to $P^\CC$.
Now the forms $\theta_k$ also descend to connection 1-forms on $P^\CC$ and we can define an almost complex structure on $P^\CC$ as $I(dr_k^2)=\theta_k$ and on the horizontal co-vectors is just a pull-back of the complex structure on the base. It defines the standard  complex structure on the fibres $(\CC^*)^n$. Its integrability follows from the fact that $\omega_k$ are of type  $(1,1)$ (see \cite{GGP}). The symplectic form is $\omega= \sum d(r_k^2\theta_k) + \pi^*(\omega_M)$, where $\omega_M$ is a K\"ahler form which is positive enough to ensure that $\omega$ is non-degenerate in the horizontal directions for almost all $x_i$. Now it is clear that for a basis of vertical vector fields $X_k$ which are defined by the $T^n$ action and satisfy $\theta_i(X_j)=\delta_i^j$, the moment map is $\mu = (r_1^2,...,r_n^2)$ as a $\RR^n$-valued function on $P^\CC$. So it is clear that for $c\in \RR^n$ where all coordinates are positive, $\mu^{-1}(c) \equiv P$ where we identify $P$ with the set of points in $P^\CC$ with $r_k=1$ for all $k$. Then it is clear that $P^\CC//T^n \equiv M$.

Q.E.D.

We can see that the reduced symplectic form depends on the level $c$ and is integral whenever $c$ satisfies
some integrality condition - which will provide the quantum condition for the correspondence in the higher rank symmetric spaces.

 We identify the reduced symplectic form on $P^\CC//T^n$ in the following way:

 \begin{co}\label{symplform}
 In the notations of the Theorem \ref{toricsymplectic} and its proof, the symplectic form on $P^\CC$ is given by $\omega = d(\sum x_i^2\theta_i)+\pi^*\omega_M$, and the reduced symplectic form on $P^\CC//T^n = \mu^{-1}(c_1,...,c_n) \equiv M$ for a generic choice of $(c_1,...,c_n)$ is $\tilde{\omega}_M = \sum c_i^2d\theta_i+\omega_M$.
 \end{co}

\section{Symmetric spaces of general rank}

 Now we apply the Corollary and the Theorem of the previous Section to the space parametrizing the maximal totally geodesic tori of a Riemannian symmetric space. Let as before $M=G/K$ be a symmetric space with $G$ compact and semisimple. Every maximal totally geodesic tori is tangent to a translated maximal commutative subspace of $\mfm$. Denote again by  $\mfa$  one such fixed subspace. Also $L$ is the connected subgroup of $G$ with Lie algebra $\mfl = \mfm+\mfa$, where $\mfm$ is the centralizer of $\mfa$ in $\mfk$. We can also write $L=MA$ where $M$ and $A$ are the corresponding Lie groups (see \cite{Gi}, \cite{Go}). Then $G/L$ is a generalized flag manifold parametrizing the maximal totally geodesic tori in $M$. As such it caries a natural complex structure, which depends on the choice of a Cartan subalgebra of $\mfg^\CC$ and a partial order in it which determines a positive Weyl chamber and it defines a positive Weyl chamber in $(\mfa^\CC)^*$. The later is dual to the cone of restricted dominant weights. Then as a complex manifold $G/L$ is equivalent to $G^\CC/M^\CC A^\CC N^\CC$ and has a principle $A^\CC$-bundle $G^\CC/M^\CC N^\CC \rightarrow G^\CC/M^\CC A^\CC N^\CC$ with total space - the {\it horospherical manifold} $\Theta$. Since $A^\CC$ is the complexification of the real torus $T^r=A$ and can be identified with the cotangent bundle $T^* T^r$, then $\Theta$ can be identified with the total space of the vertical (co)tangent bundle of the principal bundle $G/M \rightarrow G/L$ with fiber $A$. In case the rank of $M$ is $r=1$, this is just $T^*M$. Since the characteristic classes of the bundle $G/M\rightarrow G/L$ are determined via transgression by the simple roots in $\mfa^*$, then we can apply the constructions of the previous Section.

 \begin{te}\label{main}
 Let $M=G/K$ be a compact Riemannian symmetric space of rank $k$  with $G$ semisimple and let $\theta_1,...,\theta_k$ be the basis of fundamental weights that is dual to the simple restricted roots of $\mfa^\CC$. Let $\Theta$ be the associated horospherical manifold and  $\Theta\rightarrow G/L$ be the corresponding principal $(\CC^*)^k$-bundle, where $G/L$ is the quantization space of $G/K$. Let  $\omega_M=\frac{\i}{2\pi} d\rho$ be the 2-form on $G/L$ representing  $\frac{1}{2}c_1(G/L)$, so $\rho$ is  the half sum of the positive roots in $\mfg^c$ vanishing on $\mfl^c$ (as in \cite{AlP, BFR}). Then there exists a symplectic form $\omega$ on $\Theta$ with the following properties:

 i) There are positive numbers $n_i$ such that the reduced form $\tilde{\omega} $  on $\Theta//T^k$ corresponding to $\omega$ via the Marsden-Weinstein reduction is  $\tilde{\omega} = \sum_{i=1}^k n_id\theta_i + \omega_{M}$,  on $G/L$.

 ii) When the 2-form $d\alpha$ for $\alpha = n_1\theta_1 + n_2\theta_2 +...+n_k\theta_k$ is integral (up to a factor of $2\pi$) and determines a dominant weight, then the corresponding quantum bundle $\mathcal{L}$ defined by $\tilde{\omega}\in c_1(\mathcal{L})$ has the property that its space of  holomorphic sections $H^0(G/L, \mathcal{O}(\mathcal{L}))$  is an irreducible unitary representation of $G$ with highest weight $\alpha$.

iii) The complexified eigenspaces of the Laplace-Beltrami operator $\Delta_M$ corresponding to the eigenvalue $\lambda_{\alpha} = ||\alpha+\rho_{\mfa}||^2 - ||\rho_{\mfa}||^2 $
 on $M$ have dimension equal to the sum over all $\alpha$ with $||\alpha+\rho_{\mfa}||^2 = \lambda_{\alpha}+||\rho||^2$ of the dimensions of $H^0(G/L, \mathcal{O}(\mathcal{L}))$ defined in ii). All eigenvalues of $\Delta_M$ are equal to $\lambda_{\alpha}$ for some $\alpha$.

 \end{te}

 {\it Proof:} Because $c_1(G/L)>0$ for generalized flag manifolds, $\Omega_M$ is positive definite and K\"ahler.The form $\omega$ is closed and  since $\theta_i$ correspond to a basis of the positive Weyl chamber in $\mfa^\CC$, we see that the $\omega$ is also K\"ahler as a sum of a positive and non-negative form.
 %since
 %form for the complex structure on $\Theta$ because it is constructed via principal bundle structure with the aid of the connection form $\theta_i$. {\it more deteils later!!}

 The reduced form coincides with of $\tilde{\omega}_M$ by Corollary \ref{symplform}, which proves $i)$.

 The quantum line bundle is well defined since for generalized flag manifolds the Picard group is isomorphic to $H^2(G/L, \ZZ)$. Then $ii)$ follows from Borel-Weil Theorem. Finally $iii)$ is valid in view of the fact that the eigenspaces of the Laplace-Beltrami operator are sums of irreducible $G$-modules.

Q.E.D.

\begin{ex} \label{ex-su-so}The space $M=SU(3)/SO(3)$.\end{ex}

The space $SU(3)/SO(3)$ has rank 2 - which is the rank of $SU(3)$ (i.e. it is of maximal rank). If $Z=X+iY \in su(3)$, then $X\in so(3)$ and the decomposition of the Lie algebra $su(2)$ is determined by $su(3)=so(3)+\mfp$ where $\mfp$ is the space of purely imaginary matrices. The space $\mfa$ is given by the diagonal ones and the Lie algebra of $M$ is trivial. In particular the generalized flag manifold $G/L$ is the standard manifold of full flags  in $\CC^3$ identified with $SU(3)/S(U(1)\times U(1)\times U(1))$. Then a choice of the simple roots is given by diagonal matrices , which up to a factor of $\sqrt{-1}$ are $\alpha_1 = (1,-1,0), \alpha_2 = (0,1,-1)$ and the other positive root is $\alpha_3=\alpha_1+\alpha_2$. So the half sum of the positive roots is $\frac{1}{2}\rho = \frac{1}{2}(2\alpha_1+2\alpha_2)= \alpha_3 = (1,0,-1)$ for two diagonal matrices $H = (h_1,h_2,h_3),H' = (h_1',h_2',h_3')$ the product $\langle H,H' \rangle= \text{Tr}(\text{ad}_H \text{ad}_{H'}) = \sum_{i<j\leq 3} (h_i-h_j)(h_i'-h_j')$.

The dominant weights are given by $k_1\alpha_1+k_2\alpha_2$ with $k_1<2k_2<4k_1$ integers. A straightforward check gives $||k_1\alpha_1+k_2\alpha_2 - \frac{1}{2}\rho||^2= 6[(k_1+1)^2-(k_1+1)(k_2+1)+(k_2+1)^2]$. In particular we see that the eigenspaces of the Laplacian on $SU(3)/SO(3)$  corresponding to $\lambda = ||k_1\alpha_1+k_2\alpha_2 - \frac{1}{2}\rho||^2-||\frac{1}{2}\rho||^2$ split into irreducible representations subspaces when the equation $x^2 - xy + y^2 = Q$ have more than one integer solution with $x,y>1, x-1<2(y-1)<4(x-1)$. Since the condition is symmetric in $x$ and $y$, every solution $(x,y)$ with $x\neq y$ will have $(y,x)$ a solution again. The number of integer solutions of this Diophantine equation is a classical number theory question. In particular when $k_1=k_2 = n^2+1$, so that $Q=n^2$  with all of the prime factors of $n$ being of the type $3k+2$, then the solution is unique and the corresponding eigenspace of the Laplacian is irreducible $SU(3)$-module. But when $Q = 8281 = 7^213^2$ for example, then the corresponding eigenspace split into a sum of 5 irreducible $SU(3)$-modules. More details and computations for $SU(3)/SO(3)$ as well as information on the other rank 2 compact symmetric spaces appears in \cite{CM}.

%The explicit calculation and description of eigenvalues of of this space $SU(3)/SO(3)$ as well as the classification of all remaining rank 2 compact symmetric spaces appears the Ph.D. %dissertation of C. Montoya  \cite{CM}.

%\textcolor{blue}{ :) it would be awesome if I could also have a citation on my thesis too! pleaseeeee, it would be extra amazing to have a citation of my thesis, even if its for a small %reason as this one. :)}

\section{Laplace eigenfunctions and holomorphic sections}

In this Section we describe a procedure to obtain an explicit algebraic expression of the eigenfunctions of the Laplace-Beltrami operator on compact symmetric spaces through harmonic polynomials.

We shortly describe the idea of the construction first. Consider the space of holomorphic sections of a line bundle in the Borel-Weil Theorem. It is identified with holomorphic functions $f$ on a principal $\CC^*$-bundle $P$ over the (generalized) flag manifold $\mathbb{F} = G/L$ such that $f(xa) = \chi(a)f(x)$, where $\chi$ is the character of the representation in $H^0(\mathbb{F}, L_{\chi})$ for the associated with $P$ line bundle $L_{\chi}$ from the Borel-Weil Theorem. The structure group of $P$ could be reduced to $S^1$ so $P$ has a structure of a cone $P\cong \RR^+\times S$, for $S$ - the total space of an $S^1$-bundle over $\mathbb{F}$. Note that $P$ is different from $\Theta$ and sometimes can be represented as its quotient. The $S^1$ action on $S$ is induced from the $\CC^*$-action on $P$ such that for $a=re^{i\theta}\in \CC^*$ we have the action $R_a(x,t) = (e^{i\theta}x,rt)$.  We note that $S$ has a Sasakian metric $g_S$ (see \cite{BG}) and there is a cone metric $g_P$ on $P$ such that $g_P = dr^2 + r^2 g_S$. Then $g_P$ is the K\"ahler cone metric - as in Boyer-Galicki approach to Sasakian geometry \cite{BG}. In particular every holomorphic function on $P$ is also harmonic. Now the relation between the Laplace-Beltrami operators on the cone $P$ and the base $S$ is $$\Delta_P(u)=\frac{\partial^2u}{\partial r^2}+ n\frac{1}{r}\frac{\partial u}{\partial r} + r^{-2} \Delta_S(u)$$ where $u=u(r,x)$ and $\Delta_S(u) $ is calculated when $S$ is embedded in $P$ as  $r=constant$. If the function $u$ is corresponding to a holomorphic section then the equivariance condition above gives for $x=e^{i\theta}y$

$$u(x,r) = u(e^{i\theta} y, r.1) = r^k e^{ik\theta}u(e^{-i\theta}x,1)$$

where $k = \chi(re^{i\theta})$. Then from the formulas we obtain $\frac{\partial u}{\partial r} = \frac{k}{r} u$ and
\begin{equation}\label{eq1}
\Delta_Su = \lambda u
\end{equation}
when $u(x) = u(x,1)$ and $\lambda$ depends on $k$. In particular $u$ determines an eigenfunction of the Laplace-Beltrami operator on $S$. Now, in many cases, we can pull-back the function to $\Theta$ and if this pull-back is $K$-invariant, then it will define a function on $G/K$. This function is an eigenfunction if the projection is a Riemannian submersion with totally geodesic fibers. To make this strategy work we have to resolve two problems. First we need to see when a pull-back to $\Theta$ is possible. Second, the metrics on $\mathbb{F}$ which will lead to such projection are not K\"ahler - they arise from the biinvariant metric on $G$. So we need a modification of this idea for non-K\"ahler metrics. We start with the second problem.

Recall that a Hermitian metric $g$ on a complex manifold $M$ with a fundamental form $\omega$ is called balanced, if $d\omega^{n-1} = 0$ where $n$ is the complex dimension of $M$. A result in \cite{GW} shows that a holomorphic function on a balanced manifold is again harmonic. To use this property we need a few Lemmas.

\begin{lm}\label{l1}
Suppose that $M$is a compact complex manifold of dimension $n$ with a balanced metric $g_M$ which has fundamental form $\omega_M$, i.e. $d\omega_M^{n-1} = 0$
Let $\pi: P\cong \RR^+ \times S \rightarrow M$ be a principal $\CC^*$-bundle with $U(1)$-connection 1-form $\theta$ on $S$ and a cone metric of the form $g=dr^2 + r^2(\theta^2 + \pi^*(g_M))$. Let $d\theta = \pi^*(\omega)$, where $\omega$ is a form of type (1,1), be the curvature of $S$ (and $P$). With respect to a natural complex structure $I$ on $P$ compatible with $g_P$, such $g_P$ is balanced iff $$(\omega-\omega_M)\wedge\omega_M^{n-1} = 0.$$

\end{lm}

{\it Proof:} The complex structure on $P$ is given by $I(dr) = r\theta, I(\theta)=d\log r$ and the pull-back of the complex structure on $M$ on the horizontal spaces $\ker(\theta)\cap \ker(dr)$. The fact that it is integrable follows from the condition that $\omega$ is (1,1) (see \cite{GGP}). Suppose that the complex dimension of $M$ is $n$, so $dim_{\CC} P = n+1$. Then the fundamental K\"ahler form of $g_P$ is given by $\omega_P = rdr\wedge \theta + r^2\pi^*(\omega_M)$. For convenience we write $\omega_M$ for $\pi^*(\omega_M)$ where it is not confusing.  Now $\omega_P^n = nr^{2n+1}dr\wedge\theta\wedge\omega_M^{n-1} + \omega_M^n$ and $$d\omega_P^n = nr^{2n+1}dr\wedge(\omega_M-\omega)\wedge\omega_M^{n-1}$$ since $d\omega_M^{n-1} = 0$.
{\it Q.E.D.}

%{\bf Check for a factor of 2 in $\omega_M$ above!}

%{\it May be: Now the whole construction shows that by dimensionn count, all eigenfunctions arise in this manner.}

Then we have

\begin{lm}\label{l2}
Let $P \cong \RR^+\times S$ be a principal $\CC^*$-bundle over a generalized flag manifold $M = G/L$ where $G$ is a compact simply-connected and semisimple Lie group and $L$ is a centralizer of a torus. Assume that $S$ has a connection 1-form $\theta$ which is $G$-invariant and its curvature $\omega$  has a cohomology class $[\omega]$ such that $[\omega]/2\pi \in H^2(M,\ZZ)$ and is not an integer multiple of another class. Then there is a projection $\pi_1: G\rightarrow S$ which is a factor-bundle. With respect to the natural complex structure,  $P$ admits a balanced cone metric $g_P = dr^2+r^2g_S$ with an  induced $g_S$ metric on $S$, such that the projection $\pi$ is a Riemannian submersion with totally geodesic fibers when $G$ is equipped with its biinvariant metric, after possible rescaling.

\end{lm}

{\it Proof:} The fact that there is such a projection follows from \cite{G2}. For every invariant $g_M$, the Hodge-dual $*d(\omega_M)^{n-1}$ of $d\omega_M^{n-1}$ is an invariant 1-form and the Euler characteristic of $M$ is positive, so $g_M$ is balanced. From Lemma \ref{l1} we see that both  $\omega_M^{n}$  and $(\omega\wedge\omega_M^{n-1})$ are proportional to the invariant volume form on $M$. Which means that up to a rescaling of the metric on $M$ we could make them equal, so $g_P$ is balanced.

\begin{lm}\label{l3}
Every holomorphic function on $P$ is harmonic with respect to the metric $g_P$ from Lemma \ref{l2}.
\end{lm}

{\it Proof:} The result follows for example from \cite{GW}.

% Proposition 4.2
Note that most of the results in the Lemmas above are valid for non-integrable almost complex structures. However, we are focusing on the integrable case, since we need the conditions for the Borel-Weil Theorem to be satisfied in order to provide the relation to the Laplace-Beltrami eigenfunctions on $G/K$

\begin{lm}\label{l4}
A harmonic function $F$ on $P$ which satisfies $f(x,r)=r^kf(x,1)$ induces an eigenfunction of the Laplace-Beltrami operator on $S$.
\end{lm}

{\it Proof:} It follows from (\ref{eq1}) and the calculations there.

Now we consider the problem of existence of a pull-back of a function to $\Theta$. Denote by $L_{ss}$ the subgroup of $G$ with Lie algebra $\mfl_{ss} = [\mfl,\mfl] = [\mfm.\mfm]$ which is the semisimple part of $\mfl$. Consider $G/L_{ss}$ as a $T^k$-principal bundle over the flag manifold $G/L$. It has characteristic classes given by $ \gamma_i = \frac{1}{2\pi}d\omega_i $ where $\omega_i, i=1,2..k$ are the fundamental weights. Any principal $S^1$-bundle can be characterized topologically by its first Chern class, which is a positive integer combination of these. Let $S$ be determined by $c_1(S)= \sum n_i\gamma_i$, where $n_i$ are positive integers. According to Lemma 3 in \cite{G2}, if $\text{gcd}(n_1,...,n_k) = 1$, then we can find a basis of generators $\beta_1=c_1(S), \beta_2,...,\beta_k$ of $H^2(\mathbb{F}, Z)$ and they will define an equivalent principal bundle to $G\rightarrow \mathbb{F}=G/T^k$. In particular, there is a principal $T^{k-1}$-bundle $G/L_{ss} \rightarrow S$ and we can use the construction above. If this condition does not hold, then $c_1(S) = m\beta$, for some $\beta$ and $m$ positive integer, which satisfies it. Now we can replace $S$ with another bundle $\overline{S}$ with characteristic class $\beta$. By  a standard argument (for example comparing the Euler classes  - see e.g. \cite{Co}  Ex 3.26.)  $S = \overline{S}/{\mathbb Z}_m$ as a finite cover, so we have the projections $G/L_{ss}\rightarrow \overline{S}\rightarrow S$. Now $G/L_{ss}$ have two fibrations - over $S$ and over the symmetric space $G/K$. When we induce the metrics on $G/K,$ $G/L_{ss}$, and $S$ from the biinvariant metric on $G$, both fibrations are Riemannian submersions with totally geodesic fibers. For a such a Riemannian submersion $\pi:M\rightarrow N$ the relation between the Laplace-Beltrami operators on $M$ and $N$ is:

\begin{equation}\label{f1}
\Delta^M(f\circ \pi) = (\Delta^{N}f)\circ \pi
\end{equation}

for any smooth function $f$ on $N$ - see \cite{W}

%Watson JDG 1973, https://urldefense.proofpoint.com/v2/url?u=https-3A__projecteuclid.org_download_-24pdf-5F1-24_euclid.jdg_1214431482&d=DwIGaQ&c=lhMMI368wojMYNABHh1gQQ&r=mFej8Urzh8MkQKCCPimxpw&m=lZ8YOTvV9CeN9dZYn1UtUytWk3PlXv8ovKiyaEFQ9bw&s=jUqCxNlntmhBIkHwFLcjrTor57sYLFxbNK-MyQ81d6I&e= . Also Berard-Bergery and Bourguignon,  https://urldefense.proofpoint.com/v2/url?u=https-3A__projecteuclid.org_download_pdf-255C-5F1_euclid.ijm_1256046790&d=DwIGaQ&c=lhMMI368wojMYNABHh1gQQ&r=mFej8Urzh8MkQKCCPimxpw&m=lZ8YOTvV9CeN9dZYn1UtUytWk3PlXv8ovKiyaEFQ9bw&s=fDAEt-7wFDUQkR7cwKEMyYQtIrOEryB-Nd_mp8WuYlY&e= .

The above considerations lead to:

\begin{te}\label{rel}
Suppose that $G/K$ is a compact Riemannian symmetric space and $\mathbb{F} = G/L$ is the associated generalized flag manifold - the quantization space. Let $f\in H^0(\mathbb{F}, L_{\chi})$ be a holomorphic section of a positive line bundle $L_{\chi}$ over a flag manifold $\mathbb{F}$ which is considered as function on the corresponding principle bundle $P = S\times \RR^+$ with $f(za) = \chi(a)f(z)$ for $a\in \CC^*$. Let $\pi: G/L_{ss} \rightarrow S$ and $\pi_1 :G/L_{ss} \rightarrow G/K$ be the natural projections. If the pull-back of $f$ to $G/L_{ss}\times \RR^+$ via $\pi$ is $K$-invariant, then $f$ satisfies the conditions of Lemma \ref{l4} and the function $u(x) = f(x,1)$ on $S$ defines an eigenfunction $\overline{u}$ of the Laplace-Beltrami operator on the Riemannian symmetric space $G/K$  with $\pi_1^*(\overline{u})= \pi^*(u)$.
\end{te}

{\it Proof} From the Lemmas above, $u$ is an egienfunction on $S$. By the property (\ref{f1}) $\pi^*(u)$ is an eigenfunction on $G$, and by the $K$-invariance it is a pull-back of a function on $G/K$. Then again by (\ref{f1}) $\overline{u}$ is an eigenfunction on $G/K$.
{\it Q.E.D.}

\begin{re}
Sometimes the conditions of the Theorem are too strong to define all eigenfunctions. Using the Cartan embedding $i: G/K\rightarrow G$, we see that we may use functions on $G$ depending only on the parameters defining the image of $G/K$. We are going to use this in some of the examples below, even though it is not a general statement.  Moreover, the generalized flag manifold $\mathbb{F}$ has more than one invariant complex structures - see \cite{BH} for example. Each of them defines a set of eigenfunctions described in the Theorem. We focus on examples bin which this modified process in fact generates all of the eigenfunctions. It is likely that a similar procedure could be found for most of the irreducible compact symmetric spaces.

\end{re}

We mention briefly the relation of the construction to the so-called spherical representations.  A representation $\pi$ of group $G$ in a vector space $V$ (with respect to the  Riemannian symmetric space $G/K$) is called spherical, if $V$ contains a vector, fixed by all operators in $\pi(K)$. Any unitary spherical representation of $G$ with a unit vector ${\bf e}$ fixed by $\pi(K)$, the function $G\ni dx \rightarrow \langle {\bf e}, \pi(x){\bf e} \rangle$ is positive-definite and spherical (\cite{H2}, Therorem 3.4). A function on  a Lie group $G$ is called positive definite, if for every $x_1,..,x_n \in G$ and $\alpha_1,...,\alpha_n \in {\mathbb C}$
we have $$\sum_{i,j} \phi(x_i^{-1}x_j)\alpha_i\overline{\alpha_j} \geq 0$$
Also $\phi$ is called {\it spherical} if it is $K$-bi-invariant (left and right) and also a common eigenfunction of all left-invariant operators on $G$, which are also right $K-$invariant. The Cartan-Helgason Theorem \cite{H2} characterizes the irreducible spherical representations as the ones for which the highest weight $\lambda: \mathfrak{h} \rightarrow \mathbb{C}$ satisfies $$\lambda(i(\mathfrak{h}\cap\mathfrak{k}) = 0$$ and $$\frac{\langle \lambda,\alpha \rangle}{\langle \alpha,\alpha \rangle} \in \mathbb{Z}^+, \forall \alpha \in \Sigma^+$$

Note that the irreducible representations for compact $G$ are characterized by the second condition, with $\frac{2 \langle \lambda,\alpha \rangle}{\langle \alpha,\alpha \rangle} \in \mathbb{Z}^+$ instead of $\frac{\langle \lambda,\alpha \rangle}{\langle \alpha,\alpha \rangle} \in \mathbb{Z}^+$. In particular, when $G/K$ is a Riemannian symmetric
space of maximal rank, the first condition is trivial, so the spherical representations form``half" of all irreducible representations - we call them {\it even}.

In fact, the relation between the Borel-Weil theorem and the Laplace-Beltrami eigenfunctions can  potentially reveal more information. The spaces of holomorphic sections in the Borel-Weil theorem are irreducible representations and the eigenspaces on the symmetric space $G/K$ are not unless it is a CROSS. The irreducible spherical representations are characterized as the common eigenspaces of the invariant differential operators on $G/K$. So the eigenspaces are sums of spherical representations and we expect that the correspondence in the Theorem could be extended to the irreducible spherical representations. Another approach to that (see \cite{H3, Gi}) is through the integral geometry and variations of Radon transform. But such approach provides only expressions of the spherical functions in terms of integral formulas, which are not explicit in general.

\section{Examples of harmonic polynomials and eigenfunctions}

In this section we use the notation $(\cdot , \cdot)$ for the standard $SO(n)$-invariant bilinear form on $\mathbb C^N$, i.e. $(a,b) = \sum_{i=1}^N a_ib_i$.

The classical example we want to generalize is that of the sphere $S^n$. It is known that the spherical harmonics (Laplace-Beltrami  eigenfunctions), are restrictions of the harmonic polynomials on $\RR^{n+1}$. A similar description is known for $\CC P^n$. Before we present it we recall briefly the facts we need from the theory of invariant harmonic polynomials.

Let $V$ be a real or complex vector space and $G$ a group of linear transformations of $V$. Then $G$ acts on the ring of polynomials identified as the symmetric algebra $S(V^*)$.  If $f\in S(V^*)$ is a polynomial and $X\in V$, then the directional derivative $\partial(X)$ acts on $f$ as $$(\partial(X)f)(Y) = (\frac{d}{dt}f(Y+tX))|_{t=0}$$ This extends to a map $L$  from $S(V)$ to the algebra of differential operators on $V$ which is an isomorphism. Take a positive bilinear form $B$ on $V$ and define with it the isomorphism $B:V\rightarrow V^*$. The space $S(V^*)$ has a bilinear form $\langle \langle , \rangle \rangle$ defined as $$\langle \langle p,q \rangle \rangle = (\partial(P) q)(0)$$ where $P$ is the image of $p$ under the isomorphism $L\circ B$. This coincides with he usual extension of $B$ to $S(V)$, so is a positive scalar product. We have the following property for polynomials $p,q,r$ and their corresponding differential operators $P,Q,R$::
\begin{equation}\label{dual}
\langle \langle p,qr \rangle \rangle = \langle \langle \partial(Q)p,r \rangle \rangle
\end{equation}
so the multiplication by $q$ is adjoint to the operator $\partial(Q)$.
Let $I(V^*)$ is the ideal generated by the invariant polynomials and $I^+(V^*)\subset I(V^*)$ is the subset of polynomials without constant term.  For the action of $G$ denote by  $H(V^*)$ the set of $G-harmonic$ polynomials $h$, i.e. $\partial(J)(h)=0$ for every invariant differential operator $J = L\circ B (j)$, $j\in I^+(V^*)$. Assuming that $G$ is compact by \cite{H2} Ch. 3, Theorem 1.1: $$ S(V^*)=I(V^*)H(V^*)$$
and from the proof we see that $$S^k(V^*) = (I^+(V^*)S(V^*))^k +H^k(V^*)$$ is an orthogonal decomposition with respect to $\langle \langle , \rangle \rangle$. We are going to use a particular case, when $I^+(V)$ is generated by one homogeneous polynomial $p$ of degree $l$. Then the multiplication by $p$ gives an embedding $P : S^k(V^*) \rightarrow S^{k+l}(V^*)$ such that we have an identification of the quotient space $$S^{k+l}(V^*)/P(S^k(V^*)) = H^{k+l}(V^*)$$ with the harmonic polynomials which in this case are just $\text{ker}(\partial(P))$.

Now we start with a preliminary example to illustrate the correspondence between eigenfunctions and holomorphic sections in the quantization space:

\subsection{Complex projective space $\CC P^n$} \label{subsec-2-cpn}
%\vspace{.1in}

The representation of the eigenfunctions in terms of harmonic polynomials is well-known in this case. Since we will use similar reasoning later, we outline the results. We  follow the exposition in \cite{Ta}.
%Takeuch "Modern theore of harrmonic functions??"

The quantization space is the generalized flag manifold $\mathbb F = SU(n+1)/S(U(1)\times U(1)\times U(n-1))$ and the horospherical manifold corresponding to $\CC P^n$ is a bundle over the flag $SU(n+1)/S(U(1)\times U(1)\times U(n-1))$, which is embedded as a quadric (1,1)-hypersurface in $\CP n\times \CP n$ (cf. Section 3.1). This bundle should correspond to $L_{p,q}$ for some non-negative integers $p,q$. On the other hand,
$$
H^0(\mathbb{F}, \mathcal{O}(L_{p,q})) = \mathcal S^{p,q}/ \left( (z,w) \mathcal S^{p-1,q-1} \right),
$$
where  $\mathcal S^{p,q}$ denotes the space of polynomials in $z,w \in \mathbb C^{n+1}$ of homogeneous degree $(p,q)$ (in particular, $(z,w) \in \mathcal S^{1,1}$). Alternatively,
$$
\mathcal S^{p,q} = \{ F(z,w) \in S^*(\mathbb C^{2n+2})\; | \;  F(\alpha z, \beta w) = \alpha^p \beta^q F(z,w), \mbox{ for all }a,b \in \mathbb C\}.
$$
Since, the space $\mathcal S^{p,q}$ is spanned by  the polynomials $F= f(z)g(w),\hspace{1mm} \text{deg}(f) = p,\hspace{1mm} \text{deg}(g) = q,$, then $H^0(\mathbb{F}, \mathcal{O}(L_{k,k}))$ is spanned by (the restriction to $\mathbb{F}$ of) the polynomials $p(z,w) = (a,z)^k(b,w)^k$ for all $a$, $b$ with $(a,b)=0$. Alternatively, by the duality formula (\ref{dual}), the space of sections of $L_{k,k}, k>1$ can be identified with the polynomials $p(z,w)$ for which$ \sum\frac{\partial^2}{\partial z_i \partial w_i}(p)=0$.
%On the other side - under the Cartan embedding a matrix in $SU(n+1)$ with the first row $z$ is mapped into a matrix with entries $z_i\overline{z_j}$ for $i\neq j$ and $1+|z_i|^2$ on the diagonal.

 On the other side, by Corollary after Theorem 14.4  in \cite{Ta},  the functions $(a,z)^k(b,\overline{z})^k$ span the harmonic polynomials on $\CC^{n+1}$ which induce the eigenfunctions for the $k$th eigenvalue $\lambda_k = 4k(n+k)$ of the Laplace-Beltrami operator on $\CC P^n$ relative to the Fubini-Study metric (cf. Theorem \ref{rk1}). Therefore, the  $\lambda_k$-eigenspace equals the span of the restriction given by $w=\overline{z}$  of all $p(z,w) = (a,z)^k(b,w)^k$ with $(a,b)=0$.

\subsection{Quaternionic projective space $\HP n$}

Now we consider the functions of $z=(z_0,...,z_n)$ and $w=(w_0,...,w_n)$ such that $q_i = z_i + w_i j$ are the quaternionic coordinates of $\HH^{n+1}$. The harmonic polynomials on $\CC^{n+1}$ which induce the eigenfunctions on $\CP n$ are precisely the ones on $\RR^{2n+2}$ which are invariant under $S^1$. And the invariant functions which generate these polynomials are exactly $z_i\overline{z_j}$, which also fits the Cartan embedding interpretation. Similarly for $\HH^{n+1} = \RR^{4n+4} = \CC^{2n+2}$, the harmonic polynomials which determine the eigenfunctions on $\HH P^n$ are the ones right-invariant under $Sp(1)=SU(2)$, and they are functions of the variables $q_l\overline{q_k} = z_l\overline{z_k} + w_l\overline{w_k} +  (z_kw_l-z_lw_k) j$. In particular, we see that the following functions, together with their conjugates, span the $Sp(1)$ invariant quadratic harmonic polynomials:

$$(a,z)(b,\overline{z})+(a,w)(b,\overline{w})$$
$$(a,z)(b,w)-(a,w)(b,z)$$
with $(a,b)=0$.

Now we relate the functions to the quantization space which is the symplectic isotropic Grassmannian $\mathbb{F}_{is}=\mathbb{F}_{is}(2, 2n+2)$. Using the construction in Section 3.2, $\mathbb{F}_{is}$ is embedded in the regular Grassmannian by a hyperplane section given by the holomorphic symplectic form  denoted by $I(z,w)$. Then  the Grassmannian is embedded in $P(\Lambda^2(\CC^{2n+2}))$ by Plucker relations. Since the Picard group of $\mathbb{F}_{is}(2,2n+2)$ is $\ZZ$,  every line bundle is of type $\mathcal{O}(k)$ for some power of the hyperplane section bundle $\mathcal{O}_{\mathbb{F}_{is}}(1)$ which is the restriction of $\mathcal{O}_{P(\Lambda^2\CC^{2n+2})}(1)$. In particular we can consider the section as a homogeneous polynomials of degree $k$ on the variables $U_iV_j-U_jV_i$ for $U,V\in \CC^{2n+2}$ which are orthogonal to $\sum U_iV_{n+i+1}-U_{n+i+1}V_i = I(U,V)$ relative to $\langle \langle \cdot , \cdot \rangle \rangle$. Alternatively, these polynomials are the ones  in the kernel of $ \Box = \sum(\frac{\partial^2}{\partial U_i\partial V_{n+i+1}} - \frac{\partial^2}{\partial V_i\partial U_{n+i+1}})$. Using similar reasoning as the one in the case of the generalized flag manifold $\mathbb{F}$ related to $\CC P^n$, we observe that all such polynomials are linear combinations of $p(U,V)=l_{AB}(U,V)^k$ where
$$l_{AB}(U,V) = (A,U)(B,V)-(B,U)(A,V)$$
for $A,B$ satisfying certain conditions that will be determined later. To describe the conditions, we first introduce some notation and conventions. With a slight abuse of notation, consider the elements in $\CC^{2n+2}$ as pairs $\{a,b \}$ of  elements $a,b \in \CC^{n+1}$. Then if $A=\{a,b\}$ and $B=\{ c, d\}$  are  in $\CC^{2n=1}$, $I(A, B) = I( \{a,b\}, \{c,d\}) ) = (a,d)- (b,c)$  and $ (A, B) = ( \{a,b\}, \{c,d\}) ) = (a,c)+(b,d)$.  If $U =\{z,w\}, V = \{u,v\}$, then

$$\Box (l_{AB}(U,V)) = (a,d) - (b,c) = I(A,B)$$

To calculate $\Box(l_{AB}(U,V)^k)$ we use $U,V, A, B$ as above and

$$  \frac{\partial^2}{\partial U_i\partial V_{n+i+1}} (l_{AB}(U,V)^k)= \frac{\partial^2}{\partial z_i\partial v_{i}}(l_{AB}(U,V)^k) =
 $$
 $$ \frac{\partial^2}{\partial z_i\partial v_i}([(a,z)+(b,w)][(c,u)+(d,v)] - [(a,u)+(b,v)][(c,z)+(d,w)])^k=$$
 $$ \frac{\partial}{\partial z_i}k(l_{AB}(U,V)^{k-1})(d_i[(a,z)+(b,w)]-b_i[(c,z)+(d,w)] = $$
 $$ k(k-1)(l_{AB}(U,V)^{k-2})(a_i[(c,u)+(d,v)]-c_i[(a,u)+(b,v)])(d_i[(a,z)+(b,w)]-b_i[(c,z)+(d.w)]) + $$
 $$k(l_{AB}(U,V)^{k-1})(a_id_i - b_ic_i) $$

Similarly
 $$
 \frac{\partial^2}{\partial V_i\partial U_{n+i+1}} l_{AB}(U,V)^k = \frac{\partial^2}{\partial w_i\partial u_i}l_{AB}(U,V)^k = $$
 $$k(k-1)(l_{AB}(U,V)^{k-2})(b_i[(c,u)-d_[(a,u)+(b,v)])(c_i[(a,z)+(b,w)]-a_i[(c,z)+(d,w)]) + $$
 $$k(l_{AB}(U,V)^{k-1}) (c_ib_i - a_id_i) $$

 Now subtracting the two identities leads to:

$$ \Box (l_{AB}(U,V)^k) = (k^2+k)(l_{AB}(U,V)^{k-1})\sum_i(a_id_i - b_ic_i) = (k^2+k)(l_{AB}(U,V)^{k-1})((a,d)-(b,c))  $$

Hence we have the following result.
\begin{lm}\label{l5}
The holomorphic sections of the quantization bundle $\mathcal{O}_{\mathbb{F}_{is}}(k)$ are the linear combinations of the polynomials $p(U,V)  = l_{AB}(U,V)^k$ satisfying
$$(a,d)-(b,c)=0$$
\end{lm}

Now for the functions which restrict to harmonic polynomials on $\HH^n$ we take $$V = j U = \{\overline{w},-\overline{z} \}$$
This transforms $\Box$ into the Laplacian of $\HH^n$,
so that $p(U,jU)$ are pull-backs of the eigenfunctions for the $k$-th eigenvalue $\lambda_k$. More precisely we have:
\begin{te} The eigenspace of the Laplace-Beltrami operator on $\HH P^n$ corresponding to the $k$th eigenvalue is the span of the functions whose pull-back to $\HH^{n+1}$ equals
$$
p(U,jU) = ([(a,z)+(b,w)][(c,\overline{w})- (d,\overline{z})] - [(c,z)+(d,w)][(a,\overline{w})-(b,\overline{z})])^k,
$$
where $(a,d)-(b,c) = 0$.
\end{te}

{\it Proof:} We use the reasoning in \S \ref{subsec-2-cpn}. It remains to prove that the span $\mathcal{H}$ of all $p(U,jU)$ equals all of the space of harmonic polynomials on $\HH^n$. For this, we  adopt the proofs of Theorem 14.2 and 14.4 in \cite{Ta}, to the case  of  $\HH^n$  and the quantization bundle $\mathcal{O}_{\mathbb{F}_{is}}(k)$. Lemma \ref{l5} provides a description of the space of  holomorphic sections of $\mathcal{O}_{\mathbb{F}_{is}}(k)$ as a span of a set of functions. This space is a simple $Sp(n+1)$-module of highest weight $k \Lambda_2$, where $\Lambda_2$ is the second fundamental weight of $Sp(n+1)$. On the other hand, the $\lambda_k$-eignespace is also a simple $Sp(n+1)$-module of highest weight $k \Lambda_2$. Therefore the two spaces are isomorphic and set defined in the statement in the theorem corresponds to the spanning set defined in Lemma \ref{l5}. {\it Q.E.D.}

% to show that $\mathcal{H}$ is the irreducible representation of $Sp(n+1)$ on the space which is the highest component of the symmetric power of the basic representation $\HH^{n+1} = \CC^{2n+2}$. \textcolor{blue}{ The details are left to the reader.}

We note that a similar set of harmonic polynomials for $\HH^{n+1}$ defining eigenfunctions on $\HH P^n$ was found in \cite{Fu1}, Proposition 3.4, but its span was not explicitly discussed there.

%In order to see that all harmonic $Sp(1)$-invariant polynomials arise from the identity (above) we compare the dimensions with the multiplicity of the corresponding eigenvalue????but how???

\subsection{The space $SU(3)/SO(3)$}

In this example we extend the results from \cite{GSS} and find a generating set for all eigenfunctions of the Laplace-Beltrami operator on $SU(3)/SO(3)$.  Following the notations there we denote by $z,w$ etc. matrices in the Lie groups $SU(n)$ or $SO(n)$ (so $z\in SU(n)$). If $z\overline{z}^T = I$ with $I$ being the identity matrix) and by $Z,W$ matrices in the corresponding Lie algebras. Denote by $z_{ij}$ the entries of $z$. The standard metric on $SU(n)$ is given by $g(Z,W) ={\rm Re}\left(\Tr (Z\overline{W}^T)\right)$. Then the Laplace-Beltrami operator $\Delta$ on $SU(n)$ (denoted by $\tau$ in \cite{GSS}) satisfies $\Delta(fg) = \Delta(f) g + k(f,g)+\Delta(g) f$  where $k(f,g) = g(\nabla f, \nabla g)$.

For $a \in {\mathbb C}^{n}$ we set
$$
\phi_a = \sum_{j,\alpha} a_{j}a_{\alpha}\Phi_{j\alpha} = \Tr(z^T a a^T z) = (z^Ta, z^Ta).
$$
By Proposition 4.1 in \cite{GSS},  $\phi_a$ is an $SO(n)$-invariant $\Delta$-eigenfunction on $SU(n)$. Using the reasoning of  \cite{GSS},  we may show that
$$\tilde{\phi}_a (z) = \sum_{j,\alpha} a_{j}a_{\alpha}\overline{\Phi_{j\alpha}}  = \Tr(\bar{z}^T a a^T \bar{z})  =  (\bar{z}^Ta, \bar{z}^Ta)
$$
 is also  an $SO(n)$-invariant $\Delta$-eigenfunction on $SU(n)$ with the same eigenvalue as $\phi_a$.

More generally, for $a, b \in {\mathbb C}^{n}$ with $(a,b) = 0$ and $p,q \geq 0$, we have that $ \phi_a^{p} \tilde{\phi}_b^q$ is also an $SO(n)$-invariant $\Delta$-eigenfunction on $SU(n)$ of the same eigenvalue. Indeed, we can show that
$$
k(\Phi_{j\alpha}, \overline{\Phi_{k\beta}}) = -2\delta_{r\alpha}\delta_{j\beta} - 2\delta_{\alpha\beta}\delta_{kj} + \frac{4}{n}\Phi_{j\alpha}\overline{\Phi_{k\beta}}
$$
which implies that $ k(\phi_a, \tilde{\phi}_b) = -4(a,b)^2 + \frac{4}{n} \phi_a\tilde{\phi}_b$. Then using that $(a,b)=0$
and the general formula $k(f^p, h^q)  = pqf^{p-1}h^{q-1}k(f,h)$ and the fact that $k(\phi_a, \tilde{\phi}_b) = 0$, we prove that
$ \phi_a^{p} \tilde{\phi}_b^q$ are eigenfunctions. Denote by $\lambda_{p,q}$ the eigenvalue of  $\phi_a^{p} \tilde{\phi}_b^q$, $(a,b)=0$. Formulas for $\lambda_{p,q}$ in the case $n=3$ are listed in Example \ref{ex-su-so}.

Consider the space $\mathcal H^{p,q} = \Span\left\{ \phi_a^{p} \tilde{\phi}_b^q \; | \; a,b \in {\mathbb C}^{n}, (a,b) = 0 \right\}$. Then  $\mathcal H^{p,q}$ has an $SU(n)$-module structure via the formula
\begin{equation} \label{action}
A \cdot ( \phi_a^{p} \tilde{\phi}_b^q) =  \phi_{\bar{A}a}^{p} \tilde{\phi}_{Ab}^q
\end{equation}
(recall that $\bar{A} = \left(A^T\right)^{-1}$). We retain the notation from \S \ref{subsec-cpn} and \S \ref{subsec-2-cpn}. In particular, we consider the flag $\mathbb F = SU(n)/S(U(1)\times  U(1) \times U(n-2) )$ to be embedded as a quadratic hypersurface in  $\CP {n-1} \times \CP {n-1}$. Also, the holomorphic line bundles over $\mathbb F$ are denoted by $L_{k_1,k_2}$ (see the proof of Theorem \ref{th-cp}). Denote by $\Lambda_i$ the $i$th fundamental weight of $SU(n)$.

\begin{te}
The spaces $\mathcal H^{p,q}$ and $H^0(\mathbb{F}, \mathcal{O}(L_{2p,2q}))$ are both isomorphic to the simple highest weight $SU(n)$-module of highest weight $2q\Lambda_1 + 2p \Lambda_{n-1}$.  In the case $n=3$, we have that the $\Delta$-eigenspace of  $SU(3)/SO(3)$ of eigenvalue $\lambda_k$ is spanned by the  functions $\phi_a^p\tilde{\phi}_b^q$ for which $(a,b)=0$ and  $\lambda_{p,q}=\lambda_k$.
\end{te}

{\it Proof:}  Recall from \S \ref{subsec-2-cpn} that
$$H^0(\mathbb{F}, \mathcal{O}(L_{2p,2q})) \simeq \mathcal S^{2p,2q} /((z,w)\mathcal S^{2p-1,2q-1} ).$$
where $\mathcal S^{2p,2q}$ is  space of polynomials in $z,w \in \mathbb C^n$ of homogeneous degree $(2p,2q)$.The space on the right hand side is a simple $SU(n)$-module of highest weight $2q\Lambda_1 + 2p \Lambda_{n-1}$, which is isomorphic to the space spanned by the functions
$h^{p,q}_{a,b}$, $a,b \in \mathbb C^n$ with  $ ( a,b)=0$, where  $$h^{p,q}_{a,b}(z,w) = ( a,z)^p (b,w)^q.$$
This  follows from  Theorem 14.4 in \cite{Ta}. Note that the discussion in  \cite{Ta} concerns polynomials in $z$ and $\bar{z}$, but since they are treated as independent variables, the same results apply for polynomials in $z,w$. The action of $A \in SU(n)$ on $h^{p,q}_{a,b}$ is similar to the one in \eqref{action}:  $A \cdot h^{p,q}_{a,b}  = h^{p,q}_{\bar{A}a,Ab} $.

After verifying that the weights of $\phi_a^{p} \tilde{\phi}_b^q $ and $h^{2p,2q}_{a,b}$ coincide, we see that the map $ \phi_a^{p} \tilde{\phi}_b^q \mapsto h^{2p,2q}_{a,b}$ leads to the desired isomorphism $\mathcal H^{p,q} \simeq H^0(\mathbb{F}, \mathcal{O}(L_{2p,2q}))$.

  Since by the Borel-Weil Theorem $h^{m,n}_{a,b}$ generate all irreducible $SU(3)$-modules, and the eigenspaces of the Laplace-Beltrami operators on Riemannian symmetric spaces are finite sums of such modules which are even by Caratan-Helgason Theorem, we obtain the result concerning $n=3$. {\it Q.E.D.}


\begin{thebibliography}{69}



\bibitem{AM} R. Abraham, J. Marsden, {\it Foundations of mechanics}, Second edition,  (with the assistance of Tudor Ratiu and Richard Cushman), Benjamin/Cummings Publishing Co., Inc., Advanced Book Program, Reading, Mass. (1978).



\bibitem{AlP} D. Alekseevsky, A.Perelomov, {\em Invariant
K\"ahler-Einstein metrics on compact homogeneous spaces},
Funct. Anal. Appl. (3) {\bf 20} (1986), 1--16.


\bibitem{BFR} M. Bordemann, M. Forger, H, R\"{o}mer, {\it Homogeneous K\"ahler Manifolds: Paving the Way Towards New Supersymmetric Sigma Models}, Comm. Math. Phys. {\bf 102} (1986), 605--647.


\bibitem{BH}  A. Borel, F. Hirzebruch, {\em Characteristic Classes and Homogeneous Spaces, I},  Amer. J. of Math. {\bf 80} (1958), 458 -- 538.

\bibitem{BG} C. Boyer, K. Galicki, {\em Sasakian geometry} Oxford Univ. Press (2009).

\bibitem{BR} F. Burstall, J. Rawnsley, {\em Twistor theory for Riemannian symmetric spaces} Lecture Notes in Mathematics v. 1424 , Springer-Verlag (1992).

\bibitem{CW}  R. Cahn, J. Wolf, {\em Zeta functions and their expansions for compact symmetric spaces of rank one},   Commentarii Math. Helvetici {\bf 5} (1976), 1--21.


\bibitem{CS} A. Cap, J. Slovak, {\em Parabolic geometries I, background and general theory}, Mathematical surveys and monographs v. 154, American Mathematical Society (2009).


\bibitem{Cl} J. Clerc,  {\em Fonctions spheriques des espaces symetriques compacts}, Trans. Amer. Math. Soc. {\bf 306} (1988), 421--431.

\bibitem{Co} E. Correa, {\em Hermitian non-K\"ahler structures on products of principal $S^1$-bundles over complex flag manifolds and applications in Hermitian geometry with torsion}, arXiv 1803.09170.


\bibitem{Cz} J. Czyz, {\em On geometric quantization and its connections with the Maslov theory}, Rep. Math. Phys. {\bf 15} (1979),  57---97.


\bibitem{FT} K. Furutani, R. Tanaka, {\em A K\"ahler structure on the punctured cotangent bundle of complex and quaternion projective spaces and its application to a geometric quantization I}, J. Math. Kyoto Univ. {\bf 34-4} (1994) 719--737.

\bibitem{Fu1} K. Furutani, {\em Quantizationn of the geodesic flow on quaternionic projective spaces}, Ann. Global Anal. Geom. {\bf 21}, (2002), 1 -- 20.

\bibitem{Fu} K. Furutani, {\em A K\"ahler structure on the punctured cotangent bundle of the Cayley projective plane.} Jean Leray '99 Conference Proceedings, 163--182,
Math. Phys. Stud., 24, Kluwer Acad. Publ., Dordrecht, 2003.

\bibitem{Gi} S. Gindikin, {\em Horospherical Cauchy-Radon transform on compact symmetric spaces}, Mosc. Math. J.  {\bf 6} (2006), 299--305,


\bibitem{Go} R. Goodman,  {\em Harmonic Analysis on Compact Symmetric Spaces: the Legacy of Elie Cartan and Hermann Weyl}, Groups and analysis, London Math. Soc. Lecture Note Ser., 354, Cambridge Univ. Press, Cambridge, (2008) 1-- 23.


\bibitem{GGP} D. Grantcharov, G. Grantcharov, Y.Poon, {\it Calabi-Yau connections with torsion on toric bundles}, J. Differential Geom. {\bf 78} (2008), no. 1, 13--32

\bibitem{G2} G. Grantcharov {\em Geometry of compact complex homogeneous spaces with vanishing first Chern class}, Adv. Math. {\bf 226} (2011), no. 4, 3136--3159.

\bibitem{GG} D. Grantcharov, G. Grantcharov, {\em Relations between Laplace spectra and geometric quantizaion of Reimannian symmetric spaces}.  J. Geom. Symm. Phys.  {\bf 51} (2019), 9--28.


\bibitem{GW} S. Gudmundson, J. Wood, {\em Harmonic morphisms between almost Hermitian manifolds} Boll. U.M.I.  B (7) {\bf 11} (1997), no. 2, suppl., 185--197.

\bibitem{GSS} S. Gudmundsson, A. Siffert, M. Sobak, {\em Explicit proper p-harmonic functions on the Riemannian symmetric spaces SU(n)/SO(n), Sp(n)/U(n), SO(2n)/U(n), SU(2n)/Sp(n)},
preprint, arxiv:2008.08555.


\bibitem{H1} S. Helgason, {\em Differential geometry, Lie groups, and symmetric spaces. Corrected reprint of the 1978 original}, Graduate Studies in Mathematics, 34, American Mathematical Society, Providence, RI (2001).

\bibitem{H2} S. Helgason, {\em Groups and geometric analysis: integral geometry, invariant differential operators, and spherical functions},   Mathematical Surveys and Monographs, 83, American Mathematical Society (1984).

\bibitem{H3} S. Helgason, {\em Geometric analysis of symmetric spaces}, 2nd ed. Mathematical surveys and monographs v.39. American Mathematical Society, Providence, RI, (2008).


\bibitem{He} H. Hess,  {\em On a geometric quantization scheme generalizing those of Kostant-Souriau and Czyz}, Differential geometric methods in mathematical physics (Proc. Internat. Conf., Tech. Univ. Clausthal, Clausthal-Zellerfeld, 1978), Lecture Notes in Phys., 139, Springer, Berlin-New York, (1981) 1--35.

\bibitem{El} E. Hristova (Gurova), {\em Geometric Quantization of $\HP n$}, M.Sc. Thesis, University of Sofia (2007).


\bibitem{MT} I. Mladenov, V. Tsanov, {\em  Geometric quantization of the geodesic flow on $S^n$}, Differential geometric methods in theoretical physics (Shumen, 1984),  World Sci. Publishing, Singapore, (1986) 17--23.


\bibitem{CM} C. Montoya, {\em Geometric quantizations related to the Laplace eigenspectra of compact Riemannian symmetric spaces via Borel-Weil-Bott Theory. } Ph.D. Dissertation (2021), Florida International University.



\bibitem{OV} A.L. Onishchik, E. Vinberg, {\em  Lie Algebras and Lie Groups III, structure of Lie groups and Lie algebras}, Enciclopedia of Mathematical Sciences c. 43, Springer-Verlag, (1994).

\bibitem{Ra}  J. Rawnsley, {\em A non-unitary pairing of polarization for the Kepler problem}, Trans. Amer. Math. Soc. {\bf 250} (1979), 167--180.

\bibitem{Serre} J. P. Serre: {\em Representations lineaires et espaces homogenes Kehleriens des groupes de Lie compacts (d'apres Armand Borel et Andre Weil)}, Seminaire Bourbaki (Paris: Soc. Math. France) 2 (100), 1995, 447--454.

\bibitem{Ta}  M. Takeuchi, {\em Modern spherical functions}, Translation of Mathematical Monographs, v. 135 American Mathematical Society (1994).

\bibitem{W} B. Watson, {\em Manifold maps commuting with the Laplacian} J. Differential Geom. {\bf 8} (1973), no. 1, 85-94

\end{thebibliography}
\end{document}